\documentclass[a4paper,10pt]{amsart}
\usepackage[arrow,matrix]{xy}
\usepackage{amsmath,amssymb,amscd,bbm,amsthm,mathrsfs,dsfont,hyperref}
\theoremstyle{plain}
\newtheorem{thm}{Theorem}[section]
\newtheorem{cor}[thm]{Corollary}
\newtheorem{lem}[thm]{Lemma}
\newtheorem{prop}[thm]{Proposition}
\newtheorem{defn}[thm]{Definition}
\newtheorem{exa}[thm]{Example}
\newtheorem{num}[thm]{}

\makeatletter

\def\s{\sum}
\def\op{\oplus}

\def\Hom{\operatorname {Hom}}
\def\RHom{\operatorname {RHom}}
\def\Ext{\operatorname {Ext}}

\def\Tor{\operatorname {Tor}}
\def\gldim{\operatorname {gldim}}
\newcommand{\darrow}[2]{%
\genfrac{}{}{0pt}{}{\raisebox{-0.9ex}{$\xrightarrow{#1}$}}{\raisebox{0.9ex}{$\xleftarrow[#2]{}$}}}

\title{\bf Koszul differential graded algebras \\ and BGG correspondence}

\author{J.-W. He}
\address{Institute of Mathematics, Fudan University, Shanghai 200433, China}
\email{hejw@fudan.edu.cn}
\author{Q.-S. Wu}
\address{Institute of Mathematics, Fudan University, Shanghai 200433, China}
\email{qswu@fudan.edu.cn}

\date{}

\begin{document}
\begin{abstract} The concept of Koszul differential graded algebra
(Koszul DG algebra) is introduced. Koszul DG algebras exist
extensively, and have nice properties similar to the classic Koszul
algebras. A DG version of the Koszul duality is proved.
When the Koszul DG algebra $A$ is AS-regular, the Ext-algebra $E$ of
$A$ is Frobenius. In this case, similar to the classical BGG
correspondence, there is an equivalence between the stable category
of finitely generated left $E$-modules, and the quotient
triangulated category of the full triangulated subcategory of the
derived category of right DG $A$-modules consisting of all compact
DG modules modulo  the full triangulated subcategory consisting of
all the right DG modules with finite dimensional cohomology. The
classical BGG correspondence can be derived from the DG version.
\end{abstract}




\maketitle

\section*{Introduction}

In his book \cite{Ma} Manin presented an open question: How to
generalize the Koszulity to differential graded (DG for short)
algebras? Attempts have been made by several authors as in \cite{PP}
and  \cite{Be}. In their terminology, a DG algebra is said to be
Koszul if the underlying graded algebra is Koszul. Koszul DG
algebras in their sense are applied to discuss configuration spaces.

In this paper, we take a different point of view. Let $k$ be a
field. A connected DG algebra over $k$ is a positively graded
$k$-algebra $A=\op_{n \ge 0}A^n$ with $A^0=k$ such that there is a
differential $d: A \to A$ of degree $1$ which is also a graded
derivation. A connected DG algebra $A$ is said to be a {\it Koszul
DG algebra} if the minimal semifree resolution of the trivial DG
module ${}_Ak$ has a semifree basis consisting of homogeneous
elements of degree zero (Definition \ref{def1}). Our definition of
Koszul DG algebra is a natural generalization of the usual Koszul
algebra. As we will see in Section 2, a connected graded algebra
regarded as a DG algebra with zero differential is a Koszul DG
algebra if and only if it is a Koszul algebra in the usual sense.
Examples of Koszul DG algebras can be found in various fields. For
example, let $M$ be a connected $n$-dimensional $C^\infty$ manifold,
and let $(\mathcal{A}^*(M)=\bigoplus_{i=0}^n\mathcal{A}^i(M),d)$ be
the de Rham complex of $M$, then $(\mathcal{A}^*(M),d)$ is a
commutative DG algebra and by de Rham theorem (\cite{Mo}) the 0-$th$
cohomology group $H^0(\mathcal{A}^*(M))\cong \mathds{R}$. Hence the
DG algebra $\mathcal{A}^*(M)$ has a minimal model $A$ (\cite{KM}) or
Sullivan model (\cite{FHT}), which is certainly a connected DG
algebra. If the manifold $M$ has some further properties (e.g.,
$M=T^n$ the $n$-dimensional torus), then the de Rham cohomology
algebra $H(\mathcal{A}^*(M))$ is a Koszul algebra. Hence the
cohomology algebra of its minimal model (or Sullivan model) $A$ is
Koszul as $A$ is quasi-isomorphic to $\mathcal{A}^*(M)$.  Then $A$
is a Koszul DG algebra by Proposition \ref{prop2}. More examples of
Koszul DG algebra will be given in Section 2. In fact, we will see
that any Koszul algebra can be viewed as the cohomology algebra of
some Koszul DG algebra.

Bernstein-Gelfand-Gelfand in \cite{BGG} established an equivalence
between the stable category of finitely generated graded modules
over the exterior algebra $\bigwedge V$ with $V=kx_0\op
kx_1\op\cdots\op kx_n$, and the bounded derived category of coherent
sheaves on the projective space $\mathbb{P}^n$. This equivalence is
now called the {\it BGG correspondence}. BGG correspondence has been
generalized to noncommutative projective geometry by several
authors. Let $R$ be a (noncommutative) Koszul algebra. If $R$ is
AS-regular, J{\o}rgensen proved in \cite{Jo} that there is  an
equivalence  between the stable category over the graded Frobenius
algebra $E(R)=\Ext_R^*(k,k)$ and the derived category of the
noncommutative analogue QGr$(R)$ of the  quasi-coherent sheaves over
$R$; Mart\'{i}nez Villa-Saor\'{\i}n proved in \cite{MS} that the
stable category of the finite dimensional modules over $E(R)$ is
equivalent to the bounded derived category of the noncommutative
analogue qgr$R$ of the coherent sheaves over $R$. Mori in \cite{Mor}
proved a similar version under a more general condition.
One of our purposes in this paper is to establish a DG version of
the BGG correspondence.
In some special case, the DG version of the BGG
correspondence coincides with the classical one as established in
\cite{BGG} and \cite{MS}.

The paper is organized as follows.

In Section 1, we give some preliminaries and fix some notations for
the paper.

In Section 2, we first propose a definition for Koszul DG algebras
(Definition \ref{def1}), then give some examples and discuss some
basic properties of Koszul DG algebras. For any connected DG algebra
$A$, we prove that if the cohomology algebra $H(A)$ is Koszul in the
usual sense, then $A$ is a Koszul DG algebra (Proposition
\ref{prop2}). The converse is not true in general.

In Section 3, we discuss the structure of the Ext-algebras of Koszul
DG algebras. For any Koszul DG algebra $A$, we prove that the
Ext-algebra $E=\Ext_A^*({}_Ak,{}_Ak)$ of $A$ is an augmented,
filtered algebra. Moreover, if $H(A)$ is a Koszul algebra, then the
associated graded algebra $gr(E)$ is isomorphic to the dual Koszul
algebra $(H(A))^!$ (Theorem \ref{thm1}). If further, ${}_Ak$ is
compact, then $E$ is a finite dimensional local algebra; when $H(A)$
is Koszul, the filtration on $E$ is exactly the Jacobson radical
filtration (Theorem \ref{thm2}). Using bar and cobar constructions,
we prove the following version of the Koszul duality on the
Ext-algebras (Theorem \ref{thm3}):

\noindent {\bf Theorem }[Koszul Duality on Ext-algebra]. Let $A$ be
a Koszul DG algebra and $E$ be its Ext-algebra. If ${}_Ak$ is
compact, then $\Ext^*_E({}_Ek,{}_Ek)\cong H(A)$.

As a corollary, we show that the Ext-algebra of a Koszul DG algebra
$A$ with ${}_Ak$ compact is strongly quasi-Koszul (\cite{GM}) if and
only if its cohomology algebra $H(A)$ is a Koszul algebra.

In Section 4, by using Lef\`{e}vre-Hasegawa's theorem in
\cite[Ch.2]{Le} (see Theorem \ref{thm4}), we establish a DG version
of Koszul equivalence and  duality (Theorems \ref{thm5} and
\ref{thm6}).

\noindent {\bf Theorem }[Koszul equivalence and duality]. Let $A$ be
a Koszul DG algebra and $E$ be its Ext-algebra. Suppose ${}_Ak$ is
compact. Then there is an equivalence of triangulated categories
between $\mathcal{D}^+(E)$ and $\mathcal{D}_{dg}^+(A^{op});$ and
there is a duality of triangulated categories between
$\mathcal{D}^b(\textrm{mod-}E^{op})$ and $ \mathcal{D}^c(A^{op})$.

Here  $\mathcal{D}^+(E)$ is the the derived category of bounded
below  cochain complexes of left $E$-modules;
$\mathcal{D}^b(\text{\rm mod-}E)$ (resp. $\mathcal{D}^b(\text{\rm
mod-}E^{op})$) is the bounded derived category of finitely generated
left (resp. right) $E$-modules; $\mathcal{D}_{dg}(A^{op})$ (resp.
$\mathcal{D}_{dg}^+(A^{op})$) is the derived category of right DG
$A$-modules (resp.  bounded below right DG $A$-modules), and $
\mathcal{D}^c(A^{op})$ is the full triangulated subcategory of
$\mathcal{D}_{dg}^+(A^{op})$ consisting of all the compact objects.

As a corollary, we show that each finite dimensional local algebra
with residue field $k$ can be viewed as the Ext-algebra of some
Koszul DG algebra. As a result, we see that the cohomology algebra
of a Koszul DG algebra may not be Koszul.

In Section 5, we introduce the concept of AS-regular DG algebra.
Based on the result obtained in Section 4, we show that the
Ext-algebra of an AS-regular Koszul DG algebra is Frobenius
(Propostion \ref{prop7} and Corollary \ref{cor6}).  We then prove a
correspondence between some quotient category of the derived
category of a Koszul AS-regular DG algebra and the stable category
of its Ext-algebra, which is similar to the classical BGG
correspondence (Theorems \ref{thm7} and \ref{thm10}).

\noindent {\bf Theorem } [BGG Correspondence]. Let $A$ be a Koszul
DG AS-regular algebra with Ext-algebra $E=\Ext_A^*(k,k)$. Then there
is a duality of triangulated categories between $\overline{\text{\rm
mod-}}\,E\!^{op}$ and $\mathcal{D}^c(A\!^{op}) /
\mathcal{D}_{fd}(A\!^{op})$ and an equivalence of triangulated
categories between $\overline{\text{\rm mod-}}\,E$ and
$\mathcal{D}^c(A\!^{op}) / \mathcal{D}_{fd}(A\!^{op}).$

Here $\overline{\text{\rm mod-}}\,E\!^{op}$ (resp.
$\overline{\text{\rm mod-}}\,E$) is the stable category of finitely
generated right (resp. left) $E$-modules.
$\mathcal{D}_{fd}(A\!^{op})$ is the full triangulated  subcategory
of the derived category of right DG $A$-modules consisting of all
the DG modules with finite dimensional cohomology.

The results above are generalized to Adams connected DG algebras in
Section 6. We show that the noncommutative BGG correspondence
between the triangulated categories established in \cite{Jo} and
\cite{MS} can be deduced from the BGG correspondence on Adams
connected DG algebras (Theorem \ref{thm9}).

\section{Preliminaries}

Throughout, $k$ is a field and all algebras are $k$-algebras;
unadorned $\otimes$ means $\otimes_k$ and Hom means Hom$_k$.

By a {\it graded algebra} we mean a $\Bbb Z$-graded algebra. An {\it
augmented} graded algebra is a graded algebra $A$ with an
augmentation map $\varepsilon: A \to k$ which is a graded algebra
morphism.
A positively graded algebra $A=\bigoplus_{n\ge0}A_n$
with $A_0=k$ is called a {\it connected} graded algebra. Let $M$ and
$N$ be graded $A$-modules. $\underline{\hbox{\rm Hom}}_A(M,N)$ is
the set of all graded $A$-module morphisms. If $L$ is a graded
vector space, $L^\#=\underline{\Hom}(L,k)$ is the graded vector
space dual.

By a (cochain) {\it DG algebra} we mean a graded algebra
$A=\bigoplus_{n\in \Bbb Z}A^n$ with a differential $d: A \to A$ of
degree 1, which is also a graded derivation. An {\it augmented DG
algebra} is a DG algebra $A$ such that the underlying graded algebra
is augmented with augmentation map $\varepsilon : A \to k$
satisfying $\varepsilon \circ d=0$. $\ker \varepsilon$ is called the
augmented ideal of $A$. A {\it connected} DG algebra is a DG algebra
such that the underlying graded algebra is connected. Any graded
algebra can be viewed as a DG algebra with differential $d=0;$ in
this case it is called a DG algebra with {\it trivial} differential.

Let $(A,d_A)$ be a DG algebra. A \emph{left differential graded
module over $A$} (DG $A$-module for short) is a left graded
$A$-module $M$ with a differential $d_M: M \to M$ of degree $1$ such
that $d_M$ satisfies the graded Leibnitz rule
\begin{align*}
d_M(a\,m)=d_A\!(a)\,m + (-1)^{|a|}a\,d_M\!(m)
\end{align*}
for all graded elements $a \in A, \, m \in M$.

A right DG module over $A$ is defined similarly. We denote
$A\!^{op}$ as the opposite DG algebra of $A$, whose product is
defined as $a \cdot b = (-1)^{|a|\cdot|b|}ba$ for all graded
elements $a,b \in A$. Right DG modules over $A$ can be identified
with DG $A\!^{op}$-modules.

Dually, by a (cochain) {\it DG coalgebra} we mean a graded coalgebra
$C=\bigoplus_{n\in \Bbb Z}C^n$ with a differential $d: C \to C$ of
degree 1, which is also a graded coderivation. A {\it coaugmented DG
coalgebra} is a DG coalgebra $C$ with a graded coalgebra map
$\eta:k\longrightarrow C$, called coaugmentation map, such that
$d\circ\eta=0.$
If $C$ is a coaugmented DG coalgebra, then $C$ has a decomposition
$C=k\op \bar{C}$, where $\bar{C}$ is the kernel of the counit
$\varepsilon_C$, which is isomorphic to the cokernel $\widetilde{C}$
of $\eta$. There is a coproduct $\bar{\Delta}:\bar{C}\to
\bar{C}\otimes \bar{C}$ defined by $\bar{\Delta}(c)=\Delta (c) -1
\otimes c - c \otimes 1$,
such that $(\bar{C},\bar{\Delta})$ is a coalgebra without counit.
$\Delta$ induces a coproduct $\widetilde {\Delta}$ over
$\widetilde{C}$. $(\bar{C},\bar{\Delta})$ and
$(\widetilde{C},\widetilde{\Delta})$ are isomorphic as coalgebras. A
coaugmented DG coalgebra $C$ is {\it cocomplete} if, for any
homogeneous element $x\in \bar{C}$, there is an integer $n$ such
that $\bar{\Delta}^n(x)=(\bar{\Delta}\otimes 1^{\otimes
n-1})\circ\cdots\circ(\bar{\Delta}\otimes 1)\circ\bar{\Delta}(x)=0$.
A right DG $C$-comodule  $N$ is a graded right $C$-comodule with a
graded coderivation $d_N$ (i.e. $\rho_N d_N = (d_N \otimes 1 + 1
\otimes d_C) \rho_N$) of degree $1$.  A cocomplete right DG
$C$-comodule is defined similarly (\cite{Le}).

For the standard facts about DG modules, semifree modules and
semifree resolutions of DG modules, etc, refer to  \cite{AFH} and
\cite{FHT}. A DG $A$-module $M$ is said to be {\it bounded below} if
$M^n=0$ for $n\ll 0$. Let $A$ be a DG algebra, $M$ and $N$ be left
DG $A$-modules, $W$ be a right DG $A$-module. Following \cite{KM}
and \cite{We}, the differential Ext and Tor are defined as
$$\text{Ext}^n_A(M,N)=H^n(\text{RHom}_A(M,N))\quad \text{and}\quad \text{Tor}^n_A(W,M)=H^n(W\otimes^L_AM)$$
for all $n\in\Bbb Z$.

Let $s$ be the suspension map (shifting map) with $(sX)^n=X^{n-1}$
for any cochain complex $X$. Thus $s^i: X \to s^iX$ is of degree $i$
for any $i \in \mathbb{Z}$.
\begin{num}{\bf Bar constructions.}\end{num}
Let $A$ be an augmented DG algebra with differential $d$. Let
$I(A)=\cdots\op A^{-1}\op\bar{A}^0\op A^1\op\cdots$ be its augmented
ideal.  Let $$
\begin{array}{ccl} B(A)&=&T(s^{-1}(I(A)))\\ &=&k\op s^{-1}(I(A))\op
s^{-1}(I(A))\otimes s^{-1}(I(A))\op [s^{-1}(I(A))]^{\otimes
3}\op\cdots.\end{array}$$ The homogeneous element $s^{-1}a_1 \otimes
s^{-1}a_2 \otimes \cdots \otimes s^{-1}a_n$ of $B(A)$ is written as
$[a_1|a_2|\cdots|a_n]$ for homogenous elements $a_1,\cdots,a_n\in
I(A)$. The coproduct $$\Delta:B(A)\to B(A)\otimes B(A)$$ is defined
by
$$\begin{array}{ccl}\Delta([a_1|a_2|\cdots|a_n])&=
&1 \otimes [a_1|a_2|\cdots|a_n]+[a_1|a_2|\cdots|a_n]\otimes1\\
&&+\sum_{1\leq i\leq n-1}[a_1|\cdots|a_i] \otimes
[a_{i+1}|\cdots|a_n],\end{array}$$ and define a counit
$\varepsilon:B(A)\to k$ by $\varepsilon|_k=1_k$ and
$\varepsilon([a_1|\cdots|a_n])=0$ for $n\ge1$. It is easy to check
that $(B(A), \Delta, \varepsilon)$ is a coaugmented graded
coalgebra.

Define $\delta_0:B(A)\to B(A)$ by
$$\delta_0([a_1|\cdots|a_n])= - \sum_{i=1}^n(-1)^{\omega_i}[a_1|\cdots|d(a_i)|\cdots|a_n],$$
and define $\delta_1:B(A)\to B(A)$ by $$\delta_1([a_1])=0 \quad
\text{and} \quad
\delta_1([a_1|\cdots|a_n])=\sum_{i=2}^n(-1)^{\omega_i}[a_1|\cdots|a_{i-1}a_i|\cdots|a_n],$$
where $\omega_i=\sum_{j<i}(|a_j|-1)$.

It is easy to see that $\delta_0^2 = \delta_1\delta_0 + \delta_0
\delta_1=\delta_1^2=0$. Set $\delta=\delta_0+\delta_1$.  Then
$\delta$ is a differential and $(B(A),\delta)$ is a coaugmented DG
coalgebra, which is called the {\it bar construction} of $A$.

Let $(M, d_M)$ be a right DG $A$-module. The {\it bar construction
of $M$} is the complex $B(M;A)=M \otimes B(A)$ with differential
$\delta=\delta_0+\delta_1$, where
$$\begin{array}{ccl}\delta_0(m[a_1|\cdots|a_n])&=&d_M(m)[a_1|\cdots|a_n]\\
&&-\sum_{i=1}^n(-1)^{\omega_i+|m|}m[a_1|\cdots|d_A(a_i)|\cdots|a_n],\end{array}$$
and
$$\begin{array}{rcl}\delta_1(m)&=&0,\\
\delta_1(m[a_1|\cdots|a_n])
&=&(-1)^{|m|}ma_1[a_2|\cdots|a_n]\\
&&+\sum_{i=2}^n(-1)^{\omega_i+|m|}m[a_1|\cdots|a_{i-1}a_i|\cdots|a_n].\end{array}$$

$B(M;A)$ is a right DG $B(A)$-comodule.
\begin{num}
{\bf Cobar constructions.}\end{num} Let $C$ be a coaugmented DG
coalgebra with differential $d$, and let $\bar{C}=\cdots\op
C^{-1}\op \bar{C}^0\op C^1\op\cdots$ be the cokernel of the
coaugmentation map. Let
$$\Omega(C)=T(s\bar{C})=k\op s\bar{C}\op s\bar{C}\otimes s\bar{C}\op [s\bar{C}]^{\otimes 3}\op\cdots$$ be the tensor algebra,
which is augmented.
Define $\partial_0:\Omega(C)\to\Omega(C)$ by
$$\partial_0([x_1|\cdots|x_n])=- \sum_{i=1}^n(-1)^{\kappa_i}[x_1|\cdots|d(x_i)|\cdots|x_n],$$
and $\partial_1:\Omega(C)\to\Omega(C)$ by
$$\partial_1([x_1|\cdots|x_n])=\sum_{i=1}^n\sum_{(x_i)}(-1)^{\kappa_i+|x_{i(1)}|+1}[x_1|\cdots|x_{i(1)}|x_{i(2)}|\cdots|x_n],$$
where $\displaystyle \kappa_i=\sum_{j<i}(|x_j|+1)$ and
$\displaystyle\sum_{(x_i)}x_{i(1)}\otimes
x_{i(2)}=\bar{\Delta}(x_i)$. Set $\partial=\partial_0+\partial_1$.
Then $(\Omega(C),\partial)$ is an augmented DG algebra, called the
{\it cobar construction} of $C$.

Let $(M,\rho,d_M)$ be a right DG $C$-comodule. Then we have a
composition $$\bar{\rho}:M \overset{\rho}\longrightarrow M\otimes
C\longrightarrow M\otimes \bar{C}.$$ The {\it cobar construction of
$M$} is the complex $\Omega(M;C)=M\otimes \Omega(C)$ with
differential $\partial =\partial_0+\partial_1$, where
$$\begin{array}{ccl}\partial_0(m[x_1|\cdots|x_n])&=&d_M(m)[x_1|\cdots|x_n]\\
&&\displaystyle -
\sum_{i=1}^n(-1)^{\kappa_i+|m|}m[x_1|\cdots|d_C(x_i)|\cdots|x_n],\end{array}$$
and
$$\begin{array}{ccl}\partial_1(m[x_1|\cdots|x_n])&=&\displaystyle\sum_{(m)}(-1)^{|m_{(0)}|}m_{(0)}[m_{(1)}|x_1|\cdots|x_n]\\
&&\displaystyle+\sum_{i=1}^n\sum_{(x_i)}(-1)^{\kappa_i+|m|+|x_{i(1)}|+1}m[x_1|\cdots|x_{i(1)}|x_{i(2)}|\cdots|x_n],\end{array}$$
where $\displaystyle\sum_{(m)}m_{(0)}\otimes m_{(1)}=\bar{\rho}(m)$.

$\Omega(M;C)$ is a right DG $\Omega(C)$-module.

\begin{lem} \cite[Ex. 2, P. 272]{FHT}\label{lem0} Let $A$ be an augmented DG algebra. Then there is a
quasi-isomorphism of DG algebras $\zeta:\Omega B(A)\longrightarrow
A$.
\end{lem}

\begin{lem}\label{lem2}\cite[Proposition 19.2]{FHT} The augmentation map
$$B(A;A)=A\otimes B(A)\overset{\epsilon\otimes\epsilon}
\longrightarrow {}_Ak$$ is a quasi-isomorphism, and $A\otimes B(A)$
is a semifree resolution of ${}_Ak$.  \end{lem}

Dually, we have

\begin{lem}\label{lem7} The coaugmentation map $\eta:k\longrightarrow
\Omega(C;C)=C\otimes\Omega(C)$ is a quasi-isomorphism of left DG
$C$-comodule.
\end{lem}

Since $C$ is coaugmented, $C=k\op \bar{C}$. Let
$\phi:\Omega(C;C)\longrightarrow k$ be the natural linear projection
map. Then it is a right DG $\Omega(C)$-module morphism. Since
$\phi\circ\eta=id$ and $\eta$ is a quasi-isomorphism, it follows
that $\phi$ is a quasi-isomorphism, that is, $k_{\Omega(C)}$ and
$\Omega(C;C)=C\otimes\Omega(C)$ are quasi-isomorphic as DG
$\Omega(C)$-modules.

\begin{num}{\bf Some notations.}\end{num}
Let $A$ be an augmented DG algebra. $\mathcal{D}_{dg}(A)$ stands for
the derived category of left DG $A$-modules and
$\mathcal{D}_{dg}(A^{op})$ for the derived category of right DG
$A$-modules; $\mathcal{D}^c(A)$ (resp. $\mathcal{D}^c(A^{op})$)
stands for the full triangulated subcategory of
$\mathcal{D}_{dg}(A)$ (resp. $\mathcal{D}_{dg}(A^{op})$) consisting
of all the compact objects (\cite[Sect. 5]{Ke1}). If $A$ is a
connected DG algebra, then $\mathcal{D}^c(A)$ (resp.
$\mathcal{D}^c(A^{op})$) is equivalent to the full triangulated
subcategory $\langle {}_AA\rangle$ (resp. $\langle A_A\rangle$)
generated by the object ${}_AA$ (resp. $A_A$), that is, the smallest
full triangulated subcategory containing ${}_AA$ (resp. $A_A$) as an
object and closed under isomorphisms.

Let $E$ be an algebra. The notation $\mathcal{D}^*(E)$ ($*=+,\ -,\
b$) stands for the derived category of bounded below (resp. bounded
above, bounded) cochain complexes of left $E$-modules.
$\mathcal{D}^*(E^{op})$ stands for the right version of
$\mathcal{D}^*(E)$.

\section{Koszul DG algebras}

In this section, we give a definition of Koszul DG algebras, and
discuss some basic properties of Koszul DG algebras.

First of all we recall some classical definitions and well known
results. Let $V$ be a finite dimensional vector space, and
$T(V)=k\op V\op V^{\otimes 2}\op\cdots$ be the tensor algebra over
$V$. With the usual grading, $T(V)$ is a graded algebra. A {\it
quadratic algebra} is a quotient algebra $R=T(V)/(U)$ for some
finite dimensional vector space $V$ and some subspace $U \subseteq V
\otimes V$;
the {\it quadratic dual} $R^!$ of $R$ is defined as
$T(V^*)/(U^\bot),$ where $V^*$ is the dual vector space of $V$ and
$U^\bot\subseteq (V\otimes V)^*\cong V^*\otimes V^*$ is the
orthogonal complement of $U$. A quadratic algebra $R$ is {\it
Koszul\/} (\cite{Pr}, \cite{BGS}, \cite{Sm}) if the trivial
$R$-module ${}_Rk$ admits a free resolution
$$\cdots\to Q_n\to\cdots\to Q_1\to Q_0\to {}_Rk \to 0$$ with $Q_n$
generated in degree $n$ for all $n\ge0$. If $R$ is a Koszul algebra,
then its Yoneda Ext-algebra $\text{Ext}^*_R({}_Rk,{}_Rk)\cong R^!$
(\cite{Sm}, \cite{BGS}). For more properties about Koszul algebras,
we refer to the references \cite{BGS}, \cite{Pr} and \cite{Sm}.

Ungraded Koszul algebra was defined by Green-Mart\'{i}nez Villa
(\cite{GM}). Let $E$ be a noetherian semiperfect algebra with
Jacobson radical $J$. $E$ is called a {\it quasi-Koszul algebra} if
the quotient module $E/J$ has a minimal projective resolution
$$\cdots\longrightarrow
P_{n}\overset{\delta_{n}}{\longrightarrow}P_{n-1}
\overset{\delta_{n-1}}{\longrightarrow}\cdots\longrightarrow
P_{1}\overset{\delta_{1}}{\longrightarrow}P_0\overset{\delta_0}\longrightarrow
E/J\longrightarrow 0$$ such that
$$   \text{ker}\,\delta_n\cap
J^2P_n=J \, \text{ker}\,\delta_n \ \ \text{for all $n \ge 0$}.$$
 $E$ is called a {\it strongly quasi-Koszul algebra} if $$\text{ker}\,\delta_n\cap
J^iP_n=J^{i-1}\text{ker}\,\delta_n \  \  \text{for all $i \ge 2$ and
$n \ge 0$}.$$ More properties and applications of (strongly)
quasi-Koszul algebras may be found in \cite{GM} and  \cite{Mar}. We
point out here that if $E$ (with $E/J\cong k$) is a strongly
quasi-Koszul algebra then $gr(E)$, the associated graded algebra, is
Koszul (\cite{GM}).

Now let $A$ be a connected DG algebra, and let
$I=\bigoplus_{n\ge1}A^n$. A DG $A$-module $M$ with differential $d$
is said to be {\it minimal} if $d(M)\subseteq IM$. If $M$ is a
bounded below DG $A$-module, then $M$ has a minimal semifree
resolution (\cite{KM}, \cite{MW}). Recall that a DG $A$-module $P$
is called \emph{semifree} if there is a filtration of DG submodules
$$0\subseteq P(0)\subseteq P(1)\subseteq\cdots\subseteq P(n)\subseteq\cdots$$
such that $P =\cup_{n \ge 0} \,P(n)$ and each $P(n)/P(n-1)$ is free
on a basis of cocycles.

A graded subset $E$ of a DG $A$-module $P$ is called a
\emph{semibasis} if it is a basis of the graded module $P$ over the
graded algebra $A$ and has a decomposition $E = \bigsqcup_{n \ge 0}
E^n$ as a union of disjoint graded subsets $E^n$ such that
$$d(E^0)=0 \,\,\, \textrm {and} \,\,\,d(E^n) \subseteq \op_{e \in \bigsqcup_{i<n}\!E^i} Ae
\,\,\, \textrm{for all}\,\,\, n>0.$$ A DG $A$-module is semifree if
and only if it has a semibasis (\cite[Proposition 2.5]{AFH}).

We now give a definition of the Koszulity for DG algebras.

\begin{defn}\label{def1} {\rm A connected DG algebra $A$ is called a
{\it left Koszul} DG algebra if the trivial DG module ${}_Ak$ has a
minimal semifree resolution $\varepsilon: P\to {}_Ak$ such that the
semibasis of $P$ consists of elements of degree zero.}
\end{defn}

{\it Right Koszul} DG algebra is defined similarly. The next
proposition tells us that a connected DG algebra is left Koszul if
and only it is right Koszul.

\begin{prop}\label{prop1}
Let $A$ be a connected DG algebra. The following statements are equivalent.
\begin{itemize}
\item [(i)] $A$ is a left Koszul DG algebra;
\item [(ii)] $\Ext^n_A({}_Ak,{}_Ak)=0$ for all $n\neq0$;
\item [(iii)] $\Tor^n_A(k_A,{}_Ak)=0$ for all $n\neq0$;
\item [(iv)] $A$ is a right Koszul DG algebra.
\end{itemize}
\end{prop}

\proof Using the minimal semifree resolution of the trivial module.
\quad $\square$


Let $R$ be a connected graded algebra. Suppose that
$$ \cdots\to Q_n\to\cdots\to Q_1\to Q_0\to
{}_Rk \to 0$$ is a minimal free resolution of the trivial module
${}_Rk$. If we consider $R$ as a DG algebra with trivial
differential, and view $\cdots \to Q_n \to Q_{n-1} \to \cdots$ as a
double complex by using the sign trick, then the associated total
complex (that is, $Q_0 \oplus Q_1[-1] \oplus \cdots \oplus Q_n[-n]
\oplus \cdots$) is a minimal semifree resolution of the trivial DG
module ${}_Rk$. Therefore $R$ is a Koszul algebra in the usual sense
if and only if it is a Koszul DG algebra with trivial differential.

\begin{prop}\label{prop2} Let $A$ be a connected DG algebra. If the cohomology algebra $H(A)$
is a Koszul algebra, then $A$ is a Koszul DG algebra.
\end{prop} \proof We use the Eilenberg-Moore spectral sequence
(\cite{FHT}, \cite{KM})
$$E_2^{p,q}=\text{Tor}_{H(A)}^{p,q}(k,k) = \text{Tor}^{H(A)}_{-p}(k,k)^q \Longrightarrow
\text{Tor}^{p+q}_A(k,k),$$ where
$q$ is the grading induced by the gradings on $H(A)$ and
${}_{H(A)}k$. This is a convergent bounded below cohomology spectral
sequence. Since $H(A)$ is a Koszul algebra, $E^{p,q}_2=0$ for
$p+q\neq0$. Thus $\text{Tor}^{n}_A(k,k)=0$ for all $n\neq0$.
$\square$

Before proceeding to discuss further properties of Koszul DG
algebras, we give some examples here.

\begin{exa} Let $A$ be the graded algebra $k\langle x,y\rangle/(y^2,yx)$, where
$|x|=|y|=1$. Let $d(x)=xy$ and $d(y)=0$. Then $d$ induces a
differential $d$ over $A$ and $A$ is a DG algebra. It is not hard to
check that $H(A)=k\op ky$, which is a Koszul algebra. Hence by
Proposition \ref{prop2}, $A$ is a Koszul DG algebra.
\end{exa}

The following example shows that Koszul DG algebras with nontrivial
differentials exist extensively.

\begin{exa} Each Koszul algebra $R$ is the cohomology algebra of a certain
Koszul DG algebra with nontrivial differential. In fact, $R$ can be
viewed as a connected DG algebra with a trivial differential. Then
by Lemma \ref{lem2}, $\Omega B(R)$ is quasi-isomorphic to $R$ as DG
algebras. Hence $H(\Omega B(R))\cong H(R)\cong R$. Clearly $\Omega
B(R)$ is a connected DG algebra with a nontrivial differential. By
Proposition \ref{prop2}, $\Omega B(R)$ is a Koszul DG
algebra.\end{exa}

The converse of Proposition \ref{prop2} is not true, as we will see
at the end of Section \ref{sec}. However we have the following
proposition.

\begin{prop} Let $A$ be a Koszul DG algebra. If the global dimension $\mathrm {gldim}\,H(A)\leq2$,
then $H(A)$ is a Koszul algebra.\end{prop}

\proof Let
$$Q_\bullet : \quad 0\longrightarrow Q_2\longrightarrow Q_1\longrightarrow H(A)\longrightarrow
k\longrightarrow0$$ be a minimal free resolution of the trivial
module ${}_{H(A)}k$. It is direct to check that in this case the
Eilenberg-Moore resolution (\cite{FHT}, \cite{KM}) of the trivial DG
module ${}_Ak$ arising from $Q_\bullet$ can be chosen to be minimal.
If $A$ is Koszul, then the minimal free resolution $Q_\bullet$ must
be linear and hence $H(A)$ is Koszul. $\square$

The Koszulity of DG algebras is preserved under taking
quasi-isomorphisms.

\begin{lem}\label{lem4}\cite[Proposition 4.2]{KM} Let $A$ and $B$ be DG algebras.
If there is a quasi-isomorphism of DG algebras $f:A\longrightarrow
B$, then the restriction of $f$ induces an equivalence of
triangulated categories $f^*:\mathcal{D}(B)\longrightarrow
\mathcal{D}(A)$ with the inverse functor $B\otimes^L_A-$. The same
is true for $\mathcal{D}(B\!^{op})$ and $ \mathcal{D}(A\!^{op})$.
$\square$
 \end{lem}

\begin{prop}\label{prop4} Let $A$ and $B$ be connected DG algebras.
Suppose that there is a quasi-isomorphism of DG algebras
$f:A\longrightarrow B$. If $A$ (resp. $B$) is a Koszul DG algebra,
then so is $B$ (resp. $A$).
\end{prop}

\proof If $A$ is a Koszul DG algebra, then
$\text{Ext}^n_A({}_Ak,{}_Ak)=0$ for all $n\neq0$, that is,
$\text{Hom}_{\mathcal{D}(A)}({}_Ak,{}_Ak[n])=0$ for all $n\neq0$.
Hence
$$\text{Hom}_{\mathcal{D}(B)}({}_Bk,{}_Bk[n]) \cong
\text{Hom}_{\mathcal{D}(A)}(f^*({}_Bk),f^*({}_Bk)[n])
=\text{Hom}_{\mathcal{D}(A)}({}_Ak,{}_Ak[n])=0$$ for all $n\neq0$.
Hence $\text{Ext}^n_B({}_Bk,{}_Bk)=0$ for all $n\neq0$, and $B$ is
Koszul. $\square$

\section{The Ext-algebra of a Koszul DG algebra}

In this section, we study the structure of the Ext-algebra of a
Koszul DG algebra. We prove a version of the Koszul duality on
Ext-algebra for Koszul DG algebras.

Let $P$ be a semifree DG $A$-module with a semifree filtration
$$0\subseteq P(0)\subseteq P(1)\subseteq\cdots \subseteq P(n)\subseteq \cdots.$$
We may adjust
the semifree filtration of $P$ to get a
 {\it standard filtration} of $P$ as in the
following.

Let $E$ be a semifree basis of $P$. Then as a graded $A$-module,
$P=A \otimes kE$, where $kE=\op_{e\in E}ke$ is a graded $k$-vector
space. Set inductively,
$$
\begin{aligned}
 V_{\leq 0}&= V(0)= \{v \in kE \,|\,d(v) =0 \} \
\text{and}\ F(0)=A \otimes V(0) \subseteq P,\\
V_{\leq 1}&=\{v \in kE \,|\,d(v) \in F(0)\} \ \text{and}\ F(1)=A
\otimes V_{\leq 1} \subseteq P,\\
V_{\leq n}&=\{v \in kE\,|\, d(v) \in F(n-1)\}\ \text{and}\ F(n) =A
\otimes V_{\leq n} \subseteq P. \end{aligned}
$$

Let $V(n)$ be a subspace of $V_{\leq n}$ such that $V_{\leq
n}=V_{\leq n-1} \oplus V(n)$. Then for any $0 \neq v \in V(n)$,
$d(v) \in F(n-1)\backslash F(n-2).$

Obviously, $\cup_{n \ge 0}F(n)=P$ and $F(n)/F(n-1) \cong A \otimes
V(n)$ is a free DG module over a basis of cocycles. Hence
$$0\subseteq F(0)\subseteq
F(1)\subseteq\cdots\subseteq F(n)\subseteq\cdots$$ is a new semifree
filtration on $P$, which is  called  the {\it standard semifree
filtration} of $P$ associated to the semibasis $E$.

As we will see in next example, the standard semifree filtration
depends on the choice of the semibasis.

\begin{exa} {\rm Let $A$ be a connected DG algebra such that there is an element $a\in A^1$ with
$d_A(a) \neq 0$. Let $P= Ae_0 \oplus Ae_1$ as a graded free
$A$-module with $\deg(e_i)=i$ for $i =0,1$. Define $d(e_0)=0$ and
$d(e_1)=d_A(a)e_0.$ Then $P$ is a semifree DG $A$-module with a
semifree filtration
$${\bf P}:\quad 0\subseteq P(0) \subseteq P(1)=P$$
where $P(0)=Ae_0$ and $P(1)=Ae_0 \oplus Ae_1= Ae_0 \oplus
A(e_1-ae_0)=P.$

Then $E=\{e_0, e_1\}$ and $E'=\{e_0, e_1-ae_0\}$ are two semibasis
of the semifree DG module $P$. Associated to the semibasis $E$, the
standard filtration is the original one
$${\bf P}:\quad 0 \subseteq P(0) \subseteq P(1)=P.$$

Associated to the semibasis $E'$, the standard filtration is
$${\bf F}: \quad 0 \subseteq F(0) =Ae_0 \oplus A(e_1-ae_0) = P.$$}
\end{exa}

The main reason to introduce the standard filtration is that DG
morphism preserves the standard filtration as in the following
lemma, which is needed in the proof of Theorem \ref{thm1}.

\begin{lem}\label{lem1} Let $A$ be a connected DG algebra, $M$ and $N$ be minimal semifree DG
$A$-modules with the standard filtration $0\subseteq M(0)\subseteq
M(1)\subseteq\cdots$ and $0\subseteq N(0)\subseteq
N(1)\subseteq\cdots$ respectively. If the semibasis of $M$ and $N$
consist of elements of degree 0, then any DG module morphism $f:M\to
N$ preserves the filtration.
\end{lem}
\proof  Assume that there are graded vector spaces $U(i)$ and $W(i)$
for $i\ge0$ such that $M(i)/M(i-1)=A\otimes U(i)$ and
$N(i)/N(i-1)=A\otimes W(i)$. For any $u \in U(0)$, $f(u) \in
\oplus_{i \ge 0} W(i)$ and $d(f(u))=0$ since $f$ is a cochain map.
Let $f(u)=v_{i_0}+\cdots+v_{i_t}$ with $0\neq v_{i_j}\in W(i_j)$ for
$0 \leq j\leq t$ and $i_0<i_1<\cdots<i_t$. Suppose that $t \geq 1$.
By the definition of standard filtration of $N$, $d(v_{i_j})\in
N(i_{j}-1)$ and $d(v_{i_j}) \notin N(i_{j}-2)$. However,
$0=d(f(u))=d(v_{i_0}+\cdots+v_{i_t})=d(v_{i_0}+\cdots+v_{i_{t-1}})+d(v_{i_t}).$
It follows that $d(v_{i_t})=-d(v_{i_0}+\cdots+v_{i_{t-1}})\in
N(i_{t-1}-1)\subseteq N(i_t-2)$, a contradiction. Hence $t =0$ and
$f(u)\in W(0)$, which implies $f(M(0))\subseteq N(0)$.

Now suppose $f(M(n))\subseteq N(n)$. Let $\bar{M}=M/M(n)$ and
$\bar{N}=N/N(n)$. Then $f$ induces a DG morphism
$\bar{f}:\bar{M}\to\bar{N}$. $\bar{M}$ and $\bar{N}$ are minimal
semifree modules with standard semifree filtration
$$\bar{M}(0)=M(n+1)/M(n) \subseteq \bar{M}(1)=M(n+2)/M(n) \subseteq \cdots \quad \textrm {and}$$
$$\bar{N}(0)=N(n+1)/N(n) \subseteq \bar{N}(1)=N(n+2)/N(n)\subseteq
\cdots$$
respectively. By the previous narratives, we have
$\bar{f}(\bar{M}(0))\subseteq\bar{N}(0)$, which in turn implies
$f(M(n+1))\subseteq N(n+1)$. $\square$

\begin{thm}\label{thm1} Let $A$ be a Koszul DG algebra. Then

{\rm (i)} the Ext-algebra $E=\Ext_A^*({}_Ak,{}_Ak)$
of $A$ is an augmented algebra;

{\rm (ii)} there is a filtration $${\bf F}: \quad E=F_0\supseteq
F_1\supseteq\cdots\supseteq F_n\supseteq\cdots$$ on $E$ such that
$E$ is a filtered algebra. Moreover, if $H(A)$ is a Koszul algebra,
then the associated graded algebra $gr_{\bf F}(E)$
is isomorphic to the dual Koszul algebra $(H(A))^!$.
\end{thm}

\proof (i) and the first part of (ii) may be proved by using the bar
construction of $A$. We give a direct proof here for later use.

Let $\varepsilon: P \longrightarrow {}_Ak$ be a minimal semifree
resolution of the trivial DG module ${}_Ak$. Suppose that
$$0\subseteq P(0)\subseteq P(1)\subseteq\cdots\subseteq
P(n)\subseteq\cdots$$ is a standard semifree filtration of $P$
associated to some semibasis. We have graded vector spaces $V(0),
V(1), \cdots, V(n), \cdots$ such that $P(0)=A\otimes V(0)$ and
$P(n)/P(n-1)=A\otimes V(n)$ for all $n\ge1$. By the minimality of
$P$, it is easy to see that $V(0)=k$. Since $A$ is Koszul,
$$E=\text{Ext}_A^*({}_Ak,{}_Ak)=\text{Ext}_A^0({}_Ak,{}_Ak)=\prod_{i\ge0}V(i)^*=k\op\prod_{i\ge1}V(i)^*.$$
Define a decreasing filtration ${\bf F}$ on $E$ by $${\bf F}: \quad
F_0=E\ \ \text{and}\ \ F_n=\prod_{i\ge n}V(i)^*\quad \text{for}\
n\ge1.$$ We claim that $E$ is a filtered algebra with this
filtration.  For any $x \in F_n=\prod_{i\ge n}V(i)^*$ and $y \in
F_m=\prod_{i\ge m}V(i)^*$, we still use  $x$ to denote the
corresponding DG module morphism $x:P/P(n-1)\longrightarrow {}_Ak$,
and $y$ the corresponding DG module morphism
$y:P/P(m-1)\longrightarrow {}_Ak$. Since $P/P(n-1)$ is semifree,
there is a DG module morphism $f_x:P/P(n-1)\longrightarrow P$ such
that $\varepsilon\circ f_x = x$ (\cite[Lemma 6.5.3]{AFH}.
Let $g$ be the composition
$$P\overset{\pi}\longrightarrow P/P(n-1)\overset{f_x}\longrightarrow
P,$$ where $\pi$ is the natural projection map. By Lemma \ref{lem1},
$f_x$ preserves the filtration, hence $g(P(n-1))=0$ and
$g(P(n+i))\subseteq P(i)$ for all $i\ge0$.
Let $h$ be the composition
$$P\overset{\pi}\longrightarrow P/P(m-1)\overset{y}\longrightarrow
k.$$ By definition, the product $y\cdot x$ in the algebra $E$ is the
restriction of $h\circ g$ to $\bigoplus_{i\ge0}V(i)$. Since $h\circ
g(P(n+m-1))\subseteq h(P(m-1))=0$, it follows that $y\cdot x\in
\prod_{i\ge n+m}V(i)^*$. Hence  $E$ is a filtered algebra with
filtration $\{F_n\}$.

Define a map $\epsilon:E\longrightarrow k$ by $\epsilon|_k=id_k$ and
$\epsilon|_{F_1}=0$. Since $F_1$ is an ideal, $\epsilon$ is an
algebra morphism, hence an augmentation map. (i) is proved.

Now we prove the second part of (ii). Suppose that $H(A)$ is a
Koszul algebra. The trivial $H(A)$-module ${}_{H(A)}k$ has a linear
projective resolution
$$\cdots\longrightarrow H(A)\otimes V'(n)\overset{\delta_n}\longrightarrow\cdots\overset{\delta_2}\longrightarrow
H(A)\otimes V'(1)\overset{\delta_1}\longrightarrow H(A)\otimes
V'(0)\overset{\delta_0}\longrightarrow {}_{H(A)}k\longrightarrow0.$$
The Eilenberg-Moore resolution (\cite[Proposition 20.11]{FHT}) $P'$
of the DG module ${}_Ak$ arising from the previous resolution of
${}_{H(A)}k$ is minimal. Hence $P\cong P'$ as DG $A$-modules since
$A$ is connected and then $V(i)\cong V'(i)$ as vector spaces for all
$i\ge0$. For convenience, we identify $V(i)$ with $V'(i)$ for all
$i\ge0$ and $P$ with $P'$. By the construction of the filtration
${\bf F}$ on $E$, we get $F_n/F_{n-1}\cong V(n)^*$ for all $n\ge0$.
Hence we have
\begin{equation}\label{eq1}
    gr_{\bf F}(E)\cong\bigoplus_{n\ge0}V(n)^*\cong
\text{Ext}^*_{H(A)}(k,k)
\end{equation}
 as graded vector spaces. Pick elements
$x\in V(n)^*$ and $y\in V(m)^*$. As we know, $x$ and $y$ can be
extend to be DG module maps (also denoted by $x$ and $y$
respectively) $P/P(n-1)\overset{x}\longrightarrow {}_Ak$ and
$P/P(m-1)\overset{y}\longrightarrow {}_Ak$. As before, there are
filtration-preserving DG module morphisms
$f_x:P/P(n-1)\longrightarrow P$ and $f_y:P/P(m-1)\longrightarrow P$
such that $\varepsilon\circ f_x = x$ and $\varepsilon\circ f_y = y$.
Let $g$ be the composition of the DG module morphisms
$$g:P\overset{\pi}\longrightarrow P/P(n-1)\overset{f_x}\longrightarrow P\overset{\pi}{\longrightarrow}P/P(m-1)
\overset{y}\longrightarrow k.$$ Then the product $y \cdot x \in
V(n+m)^*$ of $x$ and $y$ in $gr_{\bf F}(E)$ is the restriction of
$g$ to $V(n+m)$. Since it is filtration-preserving, $f_x$ induces a
morphism of spectral sequences
$$E_*^{p,q}(f_x):E_*^{p,q}(P/P(n-1))\longrightarrow E_*^{p,q}(P).$$
In particular, $E_1^{p,q}(P/P(n-1))=H^{p+q}A\otimes V(n-p)$ and
$E_1^{p,q}(P)=H^{p+q}A\otimes V(-p)$ for all $p\leq0$ and $p+q\ge0$.
Now we regard $x\in V(n)^*$ and $y\in V(m)^*$ as elements in
$\text{Ext}^*_{H(A)}(k,k)$. Let
$h_{-p}=\bigoplus_{p\ge-q}E^{p,q}_1(f_x)$. Then we get a commutative
diagram {\tiny$$\xymatrix{
  \ar[r] &H(A)\otimes V(n+m)\ar[d]_{h_m}\ar[r]^{\qquad\quad\delta_{n+m}}&\cdots
  \ar[r]^{\delta_{n+2}\qquad\quad}&H(A)\otimes V(n+1)\ar[d]_{h_1}\ar[r]^{\delta_{n+1}} &H(A)\otimes V(n)\ar[d]_{h_0}\ar[dr]^{\eta_x}\ar[r]&0&   \\
   \ar[r]&H(A)\otimes V(m)\ar[d]_{\eta_y}\ar[r]^{\qquad\delta_{m}}&\cdots
  \ar[r]^{\delta_{2}\qquad}&H(A)\otimes V(1)\ar[r]^{\delta_{1}} &H(A)\otimes V(0)\ar[r]^{\qquad\delta_{0}} & k \ar[r]&0\\
  &k&&&&&           }$$}
where $\eta_x$ and $\eta_y$ are graded $H(A)$-module morphisms
induced by $x$ and $y$. To avoid the possible confusion, we
temporarily denote the Yoneda product on $\text{Ext}_{H(A)}^*(k,k)$
by $y*x$. By the definition of Yoneda product, $y*x$ is equal to the
restriction of $\eta_y\circ h_m$ to $V(n+m)$. Let
$\tau_m:\bigoplus_{i=0}^mV(i)\to V(m)$ be the projection map. For
any  $v\in V(n+m)$,
$$\eta_y\circ h_m(v)=\eta_y\left(E_1^{-(n+m),n+m}(f_x)(v)\right)
=\eta_y\circ\tau_m\circ f_x(v)=g(v).$$ Hence $y*x=y\cdot x$, that
is, the products on $gr_{\bf F}(E)$ and $\text{Ext}^*_{H(A)}(k,k)$
coincide under the isomorphism in (\ref{eq1}). Since $H(A)$ is
Koszul, $\text{Ext}^*_{H(A)}(k,k)\cong(H(A))^!$. Hence $gr_{\bf
F}(E)\cong(H(A))^!$. $\square$

Let $A$ be a connected DG algebra. When the trivial module ${}_Ak$
lies in $\mathcal{D}^c(A)$, the DG algebra $A$ usually has good
properties. The following proposition is clear.

\begin{prop}\label{prop3} Let $A$ be a connected DG algebra. If ${}_Ak\in \mathcal{D}^c(A)$ and $H(A)$ is a Koszul
algebra, then $\gldim H(A)<\infty$. $\square$
 \end{prop}

\begin{thm}\label{thm2} Let $A$ be a connected DG algebra. Suppose ${}_Ak\in
\mathcal{D}^c(A)$.

{\rm(i)} If $A$ is a Koszul DG algebra, then the Ext-algebra
$E=\Ext^0_A({}_Ak,{}_Ak)$ is a finite dimensional local algebra with
$E/J=k$, where $J$ is the Jacobson radical of $E$.

{\rm(ii)} If $H(A)$ is a Koszul algebra, then $gr(E)=(H(A))^!$,
where $gr(E)$ is the graded algebra associated with the radical
filtration of the local algebra $E$.
\end{thm}


\proof We use the notations in the proof of Theorem \ref{thm1}.

(i) Since ${}_Ak\in \mathcal{D}^c(A)$, there is an integer $m$ such
that $P(m)/P(m-1)\neq0$ and $P(i)/P(i-1)=0$ for all $i> m$. Hence
the filtration ${\bf F}: E=F_0\supseteq F_1\supseteq
F_2\supseteq\cdots$ stops at the $m$-th step. By Theorem \ref{thm1},
$E$ is a filtered algebra, hence for $x\in F_1$, $x^{m+1}=0$. Thus
$E$ is a local algebra with Jacobson radical $J=F_1$ and $E/J\cong
k$.

(ii) If $H(A)$ is a Koszul algebra, then by Proposition \ref{prop3},
$\gldim H(A)<\infty$. Assume that $\gldim H(A)=n$. Then the
filtration ${\bf F}$ stops at the $n$-th step, and  $J^{n+1}=0$. By
Theorem \ref{thm1}, $gr_{\bf F}(E)\cong (H(A))^!$. If we can show
$J^i=F_i$ for all $1\leq i\leq n$, then we are done. Since $V(j)=0$
for $j \geq n+1$, $E=k\op V(1)^*\op\cdots\op V(n)^*$ and
$F_i=\bigoplus^n_{j= i}V(j)^*.$
By Theorem \ref{thm1}, $gr_{\bf F}(E)$ is generated in degree 1, so
$(F_1)^n=F_n=V(n)^*$, that is, $J^n=F_n$. Similarly, since
$V(n)^*=J^n\subseteq J^{n-1}$, $V(n-1)^*\subseteq (F_1)^{n-1} +
V(n)^* =J^{n-1}$ and $F_{n-1}=V(n-1)^*\op V(n)^*\subseteq J^{n-1}$.
On the other hand, $J^{n-1} \subseteq V(n-1)^*\op V(n)^*$. Hence
$F_{n-1}=J^{n-1}$. An easy induction shows that $J^i=F_i$ for all
$1\leq i\leq n$. $\square$

We next prove a theorem similar to the Koszul duality for Koszul
algebras \cite{BGS}.

Let $A$ be an augmented DG algebra, and let $R=B(A)$ be its bar
construction.

\begin{lem}\cite[P. 272]{FHT}\label{lem3} The map $\varphi:R^\#\longrightarrow \text{\rm End}_A(A\otimes R)$ defined by
$$\varphi(f)(1[a_1|\cdots|a_n])=\s_{i=0}^n(-1)^{|f|\omega_i}1[a_1|\cdots|a_i]f([a_{i+1}|\cdots|a_n])$$
is a quasi-isomorphism of DG algebras, where
$\omega_i=|a_1|+\cdots+|a_i|-i$. $\square$
\end{lem}

\begin{lem}\label{lem5} Let $A$ be a Koszul DG algebra and $E$ be
its Ext-algebra. If ${}_Ak\in \mathcal{D}^c(A)$, then $E^\#$ is a
coaugmented coalgebra and there is a quasi-isomorphism of DG
algebras
$$\psi: \, \Omega(E^\#)\longrightarrow A.$$
\end{lem}

\proof Let $R=B(A)$ be the bar construction of $A$. Then $R$ is a
coaugmented DG coalgebra, and is concentrated in non-negative
degrees. The graded vector space dual $R^\#$ is an augmented DG
algebra.
It follows from Lemma \ref{lem2} and Lemma \ref{lem3} that $E\cong
H(\text{End}_A(A\otimes R))\cong H(R^\#)$. Since $A$ is Koszul, $E$
is concentrated in degree zero. The last isomorphism implies that
$H^i(R)=0$ for all $i>0$. Then there is naturally a
quasi-isomorphism of coaugmented DG coalgebras
$$Z^0(R)\longrightarrow R,$$ which induces a quasi-isomorphism of
augmented DG algebras
$$R^\#\longrightarrow (Z^0(R))^\#.$$
Therefore $E\cong H(R^\#)\cong (Z^0(R))^\#$ as augmented algebras.
Since ${}_Ak\in \mathcal{D}^c(A)$, $E$ is a finite dimensional
algebra. Hence $E^\#\cong Z^0(R)$ as coaugmented coalgebras, and
there is a quasi-isomorphism of coaugmented DG coalgebras
$$E^\#\longrightarrow R.$$
This induces a quasi-isomorphism of DG algebras
$$\xi:\Omega(E^\#)\longrightarrow \Omega(R)=\Omega B(A).$$
There is also a quasi-isomorphism of DG algebras $\zeta:\Omega
B(A)\longrightarrow A$ by Lemma \ref{lem0}. Hence the composition
\begin{equation}\label{eq2}
\Omega(E^\#)\overset{\xi}\longrightarrow\Omega
B(A)\overset{\zeta}\longrightarrow A
\end{equation} gives a quasi-isomorphism of DG
algebras $\psi=\zeta\circ\xi:\Omega(E^\#)\longrightarrow A$. The
proof is completed. $\square$

\begin{thm}[Koszul Duality on Ext-algebra]\label{thm3}  Let $A$ be a Koszul DG algebra and $E$ be
its Ext-algebra. If ${}_Ak\in \mathcal{D}^c(A)$, then
$\Ext^*_E({}_Ek,{}_Ek)\cong H(A)$.
\end{thm}
\proof By Lemma \ref{lem3}, $\Omega(E^\#)=B(E)^\#$ is
quasi-isomorphic to $\text{End}_E(E\otimes B(E))$. Hence
$$\text{Ext}_E^*({}_Ek,{}_Ek)\cong H\left(\text{End}_E(E\otimes B(E))\right)\cong
H(\Omega(E^\#)).$$
It follows from Lemma \ref{lem5} that
$\text{Ext}_E^*({}_Ek,{}_Ek)\cong H(A)$. $\square$

As an application of Theorem \ref{thm3}, we have the following two
corollaries, which establish relations between Koszul DG algebras
and (strongly) quasi-Koszul algebras.

\begin{cor}\label{cor1} Let $A$ be a Koszul DG algebra. If ${}_Ak\in \mathcal{D}^c(A)$, then the following are equivalent:

{\rm (i)} The Ext-algebra $E$ of $A$ is a quasi-Koszul algebra;

{\rm (ii)} $H(A)$ is generated in degree 1. \end{cor} \proof By
Theorem \ref{thm2}, $E$ is a finite dimensional local algebra with
the residue field $k$. The equivalence of (i) and (ii) follows from
Theorem \ref{thm3} and \cite[Theorem 4.4]{GM}. $\square$

\begin{cor}\label{cor2} Let $A$ be a Koszul DG algebra. If ${}_Ak\in \mathcal{D}^c(A)$,
then the following are equivalent:

{\rm (i)} The Ext-algebra $E$ of $A$ is a strongly quasi-Koszul
algebra;

{\rm (ii)} $H(A)$ is a Koszul algebra.
\end{cor}

\proof (i) $\Rightarrow$ (ii). By \cite[Theorem 6.1]{GM} and its
proof, $\text{Ext}_E^*({}_Ek,{}_Ek)$ is a Koszul algebra. Theorem
\ref{thm3} implies that $H(A)\cong \text{Ext}_E^*({}_Ek,{}_Ek)$ is a
Koszul algebra.

(ii) $\Rightarrow$ (i). Applying Theorem \ref{thm3} again,
$\text{Ext}_E^*({}_Ek,{}_Ek)\cong H(A)$ is a Koszul algebra. By
\cite[Theorem 9.1]{GM}, $E$ is a strongly quasi-Koszul algebra.
$\square$

\section{Koszul Duality}\label{sec}

Let $B$ be an augmented DG algebra and $C$ be a coaugmented DG
coalgebra.  Lef\`{e}vre-Hasegawa (\cite[Proposition 2.2.4.1]{Le})
established an equivalence between the derived category
$\mathcal{D}(B)$ and the so called coderived category
$\mathcal{D}(C)$ when $B$ and $C$ satisfy certain conditions. Thanks
for the result of Lef\`{e}vre-Hasegawa  we can prove a version of
Koszul Duality (\cite{BGS}) for Koszul DG algebras.

Let $(B,m_B,d_B)$ be an augmented DG algebra with an augmentation
map $\varepsilon_B: B\to k$, and $(C,\Delta,d_C)$ be a coaugmented
DG coalgebra with a coaugmentation map $\eta_C:k\to C$. A graded
linear map $\tau:C\to B$ of degree 1 is called a {\it twisting
cochain} from $C$ to $B$ (\cite{HMS}, \cite{Le}) if
$$\begin{array}{c}
\varepsilon_B\circ\tau\circ\eta_C=0,\ \text{and}\\
 m_B \circ(\tau\otimes \tau)\circ\Delta+d_B\circ\tau+\tau\circ d_C=0.\end{array}$$ Let
$\Omega(C)$ be the cobar construction of $C$. The twisting cochains
from $C$ to $B$ are one to one corresponding to the DG algebra
morphisms from $\Omega(C)$ to $B$. There is a {\it canonical
twisting cochain} $\tau_0:C\to \Omega(C)$ given by $\tau_0(c)=[c]$
for any $c \in \bar{C}$ and $\tau_0(k)=0$.

Let $\tau:C\to B$ be a twisting cochain. For any right DG
$C$-comodule $N$, the {\it twisted tensor product} $N \otimes_\tau
B$ (\cite{Le}, \cite{Ke}) is the right DG $B$-module defined by

(i)  $N \otimes_\tau B= N \otimes B$ as a right graded $B$-module;

(ii) the differential $\delta = d_N \otimes 1 + 1 \otimes d_B + (1
\otimes m_B)(1 \otimes \tau \otimes 1)(\rho_N \otimes 1),$ i.e.

$$\delta(n \otimes a)=d(n)\otimes a+(-1)^{|n|}n\otimes d(a)+\s_{(n)}(-1)^{|n_{(0)}|}n_{(0)}\otimes\tau(n_{(1)})a,$$
for any homogeneous elements $n\in N$ and $a\in B$.

Dually, for any DG $B$-module $M$, the {\it twisted tensor product}
$M\otimes_\tau C$ is the right DG $C$-comodule defined by

(i)  $M\otimes_\tau C=M\otimes C$ as a vector space;

(ii) the differential $\delta = d_M \otimes 1 + 1 \otimes d_C - (m_M
\otimes 1)(1 \otimes \tau \otimes 1)(1 \otimes \Delta)$, i.e.
$$\delta(m\otimes c)=d(m)\otimes c+(-1)^{|m|}m\otimes d(c)-\s_{(c)}(-1)^{|m|}m\tau(c_{(1)})\otimes c_{(2)},$$
for any homogeneous elements $m\in M$ and $c\in C$.

Let $\text{DGmod-}B$ be the category of right DG $B$-modules and
$\text{DGcom-}C$ be the category of right DG $C$-comodules. Then
there is a pair of adjoint functors $(L,R)$ (\cite{Ke}, \cite{Le}):

$$
 \text{DGcom-}C \darrow{L=-\otimes_\tau B}{R=-\otimes_\tau C}\text{DGmod-}B.
$$

Let $C$ be a cocomplete DG coalgebra, and  DGcomc-$C$ be the
category of cocomplete right DG $C$-comodules. For any $M, N \in$
DGcomc-$C$, a DG comodule morphism $f:M \to N$ is called a {\it weak
equivalence related to $\tau$} (\cite{Ke},\cite{Le}) if $L(f): LM
\to LN$ is a quasi-isomorphism. Note that a weak equivalence related
to $\tau_0$ ($B=\Omega (C)$) is a quasi-isomorphism. But the
converse is not true in general (\cite{Ke}). Let $\mathcal{K}(C)$ be
the homotopy category of $\text{DGcomc-}C$. Equipped with the
natural exact triangles, $\mathcal{K}(C)$ is a triangulated
category. Let $\mathcal{W}$ be the class of weak equivalences in the
category $\mathcal{K}(C)$. Then $\mathcal{W}$ is a multiplicative
system. The coderived category $\mathcal{D}_{dg}(C)$ of $C$ is
defined to be $\mathcal{K}(C)[\mathcal{W}^{-1}]$, the localization
of $\mathcal{K}(C)$ at the class $\mathcal{W}$ of weak equivalences
(\cite{Ke}, \cite{Le}).  Let $\mathcal{D}_{dg}(B^{op})$ be the
derived category of right DG $B$-modules. The following theorem is
proved by Lef\`{e}vre-Hasegawa in \cite[Ch.2]{Le}, and also can be
found in \cite{Ke}.

\begin{thm}\label{thm4}Let $C$ be a cocomplete DG coalgebra, $B$ an augmented DG
algebra and $\tau:C\to B$ be a twisting cochain. Then the following
are equivalent:
\begin{itemize}
    \item [(i)] The map $\tau$ induces a quasi-isomorphism $\Omega(C)\to
    B$;
    \item [(ii)] The adjunction map $$B\otimes_\tau C\otimes_\tau B\to B$$ is a
    quasi-isomorphism;
    \item [(iii)] The functors $L$ and $R$ induce an equivalence of triangulated
    categories (also denoted by $L$ and $R$)
     $$\mathcal{D}_{dg}(C) \overset {L}{\underset {R}{\rightleftarrows}} \mathcal{D}_{dg}(B^{op}).\quad \square$$
\end{itemize}
\end{thm}

Now let $A$ be a Koszul DG algebra. Suppose ${}_Ak\in
\mathcal{D}^c(A)$. By Theorem \ref{thm2}, its Ext-algebra $E$ is a
finite dimensional local algebra with the residue field $k$. Hence
the vector space dual $E^*=E^\#$ is a coaugmented coalgebra which is
of course cocomplete. Hence all the DG $E^*$-comodules are
cocomplete. Let $C=E^*$ and $B=\Omega(C)$. Clearly, $B$ is a
connected DG algebra, and the canonical twisting cochain $\tau_0:
C\to \Omega(C)$ satisfies the condition (i) in the Theorem
\ref{thm4}. Hence we have the following equivalence of triangulated
categories
$$\mathcal{D}_{dg}(E^*)\overset {L}{\underset {R}{\rightleftarrows}}\mathcal{D}_{dg}(\Omega(E^*)^{op}).$$

Let $\mathcal{D}^+_{dg}(\Omega(E^*)^{op})$ be the derived category
of all bounded below right DG $\Omega(E^*)$-modules, that is,
consisting of objects $M$ with $M^n=0$ for $n\ll0$. Since
$\Omega(E^*)$ is connected, it is not hard to see that
$\mathcal{D}^+_{dg}(\Omega(E^*)^{op})$ is a full triangulated
subcategory of $\mathcal{D}_{dg}(\Omega(E^*)^{op})$. Similarly, let
$\mathcal{K}_{dg}^+(C)$ be the homotopy category of bounded below DG
cocomplete comodules, and let $\mathcal{D}^+_{dg}(E^*)$ be the
localization of $\mathcal{K}_{dg}^+(E^*)$ at the class of weak
equivalences $\mathcal{W}^+$ in $\mathcal{K}_{dg}^+(E^*)$
($\mathcal{W}^+$ is also a multiplicative system). One can check
that $\mathcal{D}^+_{dg}(E^*)$ is a full triangulated subcategory of
$\mathcal{D}_{dg}(E^*)$. Restricting $L$ and $R$ to the
subcategories $\mathcal{D}^+_{dg}(E^*)$ and
$\mathcal{D}^+_{dg}(\Omega(E^*)^{op})$ respectively, we get the
following proposition.

\begin{prop}\label{prop5} Let $A$ be a Koszul DG algebra and $E$ be its Ext-algebra.
If ${}_Ak\in \mathcal{D}^c(A)$, then the following is an equivalence
of triangulated categories
$$\mathcal{D}^+_{dg}(E^*)\overset {L}{\underset {R}{\rightleftarrows}}\mathcal{D}^+_{dg}(\Omega(E^*)^{op}).
\quad \square$$
\end{prop}

Since $E^*$ is concentrated in degree zero, a DG $E^*$-comodule is
exactly a cochain complex of $E^*$-comodules. Hence
$\mathcal{K}^+_{dg}(E^*)=\mathcal{K}^+(E^*)$, the homotopy category
of bounded below cochain complexes of right $E^*$-comodules. It is
not hard to see that the class $\mathcal{W^+}$ of weak equivalences
related to $\tau_0$ is exactly the class of quasi-isomorphisms.
Hence
$\mathcal{D}^+_{dg}(E^*)=\mathcal{K}^+_{dg}(E^*)[(\mathcal{W^+})^{-1}]=D^+(E^*)$,
the derived category of bounded below cochain complexes of right
$E^*$-comodules. By Proposition \ref{prop5} we have the following
proposition.

\begin{prop}\label{prop6}Let $A$ be a Koszul DG algebra and $E$ be its Ext-algebra.
If ${}_Ak\in \mathcal{D}^c(A)$, then there is an equivalence of
triangulated categories (we use the same notations of the equivalent
functors as in Prop. \ref{prop5}.)
$$D^+(E^*)\overset {L}{\underset {R}{\rightleftarrows}}\mathcal{D}^+_{dg}(\Omega(E^*)\!^{op}).
\quad \square$$
\end{prop}
Since $E$ is a finite dimensional algebra, the category of left
$E$-modules is isomorphic to the category of right $E^*$-comodules
(\cite[1.6.4]{Mon}). Hence there is an equivalence of triangulated
categories
$$\mathcal{D}^+(E)\overset {F}{\underset {G}{\rightleftarrows}} \mathcal{D}^+(E^*),$$ where
$\mathcal{D}^+(E)$ is the derived category of bounded below cochain
complexes of left $E$-modules.

By Lemma \ref{lem5}, there is a quasi-isomorphism of DG algebra
$\varphi:\Omega(E^*)\longrightarrow A$. Hence by Lemma \ref{lem4},
the following gives an equivalence of triangulated categories
$$\mathcal{D}^+(A^{op})  \darrow {\qquad\varphi^*\quad } {-\otimes_{\Omega(E^*)}^LA} \mathcal{D}^+(\Omega(E^*)^{op}).$$
Let $\Phi=(-\otimes_{\Omega(E^*)}^LA)\circ L\circ F$ and
$\Psi=G\circ R\circ\varphi^*$. We have the following theorem.

\begin{thm}[Koszul Equivalence]\label{thm5} Let $A$ be a Koszul DG algebra and $E$ be its Ext-algebra.
If ${}_Ak \in \mathcal{D}^c(A)$, then we have an equivalence of
triangulated categories
$$ \mathcal{D}^+(E) \darrow  {\Phi} {\Psi}   \mathcal{D}_{dg}^+(A^{op}).$$
\end{thm}

It is easy to see that
$\Phi({}_Ek)=L(k^{E^*})\otimes_{\Omega(E^*)}^LA=\Omega(E^*)\otimes_{\Omega(E^*)}^LA=A_A$.
Temporarily write $\langle{}_Ek\rangle$ the full triangulated
subcategory of $\mathcal{D}^+(E)$ generated by ${}_Ek$. By
restricting
 $\Phi$ and $\Psi$, we get an equivalence of triangulated
categories
$$\langle{}_Ek\rangle  \darrow {\Phi_{\text{res}}} {\Psi_{\text{res}}} \mathcal{D}^c(A).$$

\begin{lem}\label{lem6}$\langle{}_Ek\rangle= \mathcal{D}^b(\text{\rm mod-}E)$,
where $\text{\rm mod-}E$ is the category of finitely generated left
$E$-modules.
\end{lem}
\proof It suffices to show that all the finitely generated
$E$-modules are in $\langle{}_Ek\rangle$. Since $E$ is finite
dimensional, any finitely generated $E$-module is finite
dimensional. Clearly, all 1-dimensional modules are in
$\langle{}_Ek\rangle$. Let $N$ be a finite dimensional module. Since
$\rm {soc}(N) \neq 0$, we have an exact sequence
$$0\longrightarrow{}_Ek\longrightarrow N\longrightarrow
N/{}_Ek\longrightarrow0.$$ Since $\dim N/{}_Ek<\dim N$, an induction
on the dimension of $N$ implies that $N$ lies in
$\langle{}_Ek\rangle$. Hence all finitely generated $E$-modules are
in $\langle{}_Ek\rangle$. $\square$

\begin{cor}\label{cor3} Let $A$ be a Koszul DG algebra and $E$ be its Ext-algebra. If ${}_Ak\in \mathcal{D}^c(A)$, then we
have an equivalence of triangulated categories
$$\mathcal{D}^b(\text{\rm mod-}E)  \darrow {\Phi_{\text{res}}} {\Psi_{\text{res}}}
\mathcal{D}^c(A^{op}).$$
\end{cor}

Since $E$ is finite dimensional, the vector space dual $(\ )^*$
induces a duality of triangulated categories
$$\mathcal{D}(\text{\rm mod-}E) \darrow {(\ )^*} {(\ )^*} \mathcal{D}(\text{\rm
mod-}E^{op}).$$

Now, we are able to give a version of the Koszul duality for Koszul
DG algebras.

\begin{thm}[Koszul Duality]\label{thm6} Let $A$ be a Koszul DG algebra and $E$ be its Ext-algebra.
Suppose ${}_Ak \in \mathcal{D}^c(A)$. Then there is a duality of
triangulated categories
$$\mathcal{D}^b(\text{\rm mod-}E^{op})\darrow {\mathcal{F}} {\mathcal{G}}
\mathcal{D}^c(A^{op}).$$

\end{thm}

It is easy to see that
\begin{equation}\label{eq3}
\mathcal{F}(k_E)=\Phi({}_Ek)=A_A
\end{equation}and
\begin{eqnarray}
  \nonumber\mathcal{F}(E_E)& =& \Phi({}_EE^*) \\
   \nonumber &= & L((E^*)^{E^*})\otimes_{\Omega(E^*)}^LA  \\
    \nonumber & \overset{(a)}\cong & k_{\Omega(E^*)}\otimes_{\Omega(E^*)}^LA\\
  \label{eq4}  &\cong& k_{A},
\end{eqnarray}
where the isomorphism ($a$) holds, because
$L((E^*)^{E^*})=\Omega(E^*;E^*)$ which is quasi-isomorphic to
$k_{\Omega(E^*)}$ as a DG $\Omega(E^*)$-module by the narrative
below Lemma \ref{lem7}.

From the proof of above results, we  have proved in fact the
following result.

\begin{cor}\label{cor4} Let $R$ be a finite dimensional local
algebra with the residue field $k$. Then there is a duality of
triangulated categories
$$\mathcal{D}^b(R^{op})\rightleftarrows\mathcal{D}^c(\Omega(R^*)^{op});$$ and
under this duality, the trivial module $k_R$ corresponds to
$\Omega(R^*)$ and $R_R$ to $k_{\Omega(R^*)}$. $\square$
\end{cor}

The following corollary was indicated in \cite{Ke} and \cite{Le}. As
an application of Corollary \ref{cor4}, we give a proof here.

\begin{cor}\label{cor5}Let $R$ be a finite dimensional local algebra with the
residue field $k$. Then the connected DG algebra $\Omega(R^*)$ is a
Koszul DG algebra. Moreover, the Ext-algebra
$\Ext_{\Omega(R^*)}^*(k,k)$ is isomorphic to $R$.
\end{cor}

\proof By Corollary \ref{cor4}, $k_{\Omega(R^*)}$ is compact, and
$$\Ext_{\Omega(R^*)^{op}}^n(k,k)=\Hom_{\mathcal{D}^c(\Omega(R^*)^{op})}(k,k[-n])
\cong\Hom_{\mathcal{D}^b(R^{op})}(R[n],R)=0$$ if $n \neq 0$.
Therefore $\Omega(R^*)$ is a Koszul DG algebra. Moreover, the
following are algebra isomorphisms
$$\Ext_{\Omega(R^*)}^*(k,k)\cong\Ext_{\Omega(R^*)^{op}}^*(k,k)^{op}\cong
\Ext_{R^{op}}^*(R,R)\cong R.\quad \square$$

In particular, by Corollary \ref{cor5}, if $k$ is algebraically
closed, then any finite dimensional local algebra can be viewed as
the Ext-algebra of some Koszul DG algebra.

\begin{exa}\label{rem1}{\rm
Now let $V=kx\op ky\op
kz$ and $R=T(V)/T^{\ge4}(V)$. Clearly, $R$ is a finite dimensional
local algebra. Then $B=\Omega(R^*)$ is a Koszul DG algebra with
$\Ext_B^*({}_Bk, {}_Bk)=R$. Since $gr(R)\cong R$ is not a Koszul
algebra, so $R$ is not a strongly quasi-Koszul algebra. By Corollary
\ref{cor2}, the cohomology $H(B)$ can not be a Koszul algebra. Hence
the converse of Proposition \ref{prop2} is not true.}
\end{exa}

\section{BGG correspondence}

In \cite{BGG}, Bernstein-Gelfand-Gelfand established an equivalence
of categories
$$\overline{\text{grmod-}}\Lambda(V)\cong
\mathcal{D}^b(\text{Coh}\mathbb{P}^n)$$ where
$\overline{\text{grmod-}}\Lambda(V)$ is the stable category of
finitely generated graded modules over the exterior algebra
$\Lambda(V)$ of an $(n+1)$-dimensional space $V=kx_0\oplus
kx_1\oplus\cdots\oplus kx_n$, and $\mathcal{D}^b
(\text{Coh}\mathbb{P}^n)$ is the bounded derived category of
coherent sheaves over the $n$-dimensional projective space
$\mathbb{P}^n$. This equivalence is now called the {\it BGG
correspondence} in literature. A sketch of the proof of the BGG
correspondence can be found also in \cite[P.273, Ex.1]{GMa}. The BGG
correspondence has been generalized to noncommutative projective
geometry by several authors (\cite{Jo}, \cite{MS}, \cite{Mor}). Let
$R$ be a Koszul noetherian AS-Gorenstein algebra
with finite global dimension. Then its Ext-algebra $E(R)$ is a
Frobenius algebra (\cite{Sm}). A version of the noncommutative BGG
correspondence was proved in \cite{MS}, which was stated as
$$\overline{\text{grmod-}}E(R)\cong \mathcal{D}^b(\text{qgr}R^{op}),$$
where $\overline{\text{grmod-}}E(R)$ is the stable category of
finitely generated graded modules over $E(R)$ and $\text{qgr}R^{op}$
is the quotient category $\text{grmod-}R^{op}/\text{tors} R^{op}$.
Let $\mathcal{D}^b_{fd}(\text{grmod-}R^{op})$ be the
full subcategory of $\mathcal{D}^b(\text{grmod-}R^{op})$ consisting
of objects $X$ with finite dimensional cohomology groups. It is well
known that (\cite{Miy})
$$\mathcal{D}^b(\text{qgr}R^{op})=\mathcal{D}^b(\text{grmod-}R^{op})/\mathcal{D}^b_{fd}(\text{grmod-}R^{op}).$$
Hence the above BGG correspondence can be stated as
\begin{equation}\label{eq5}
\overline{\text{grmod-}}\,E(R)\cong\mathcal{D}^b(\text{grmod-}R^{op})/\mathcal{D}^b_{fd}(\text{grmod-}R^{op}).
\end{equation}
In this section, we deduce a correspondence similar to (\ref{eq5})
for AS-Gorenstein Koszul DG algebras.

First of all we recall the definition of AS-Gorenstein DG algebra.
Let $A$ be a connected DG algebra. We say that $A$ is {\it right
AS-Gorenstein} (AS stands for Artin-Schelter)  if
$\RHom_{A\!^{op}}(k,A)\cong s^nk$ for some integer $n$ (\cite{FHT1},
\cite{LPWZ},\cite{LP2}); $A$ is {\it right AS-regular} if $A$ is
right AS-Gorenstein and $k_A\in \mathcal{D}^c(A\!^{op})$. Similarly,
we define left AS-Gorenstein DG algebra and left AS-regular algebra.
We call $A$ is  AS-Gorenstein (resp., regular) if $A$ is both left
and right AS-Gorenstein (resp., regular).

\begin{prop}\label{prop9} Let $A$ be a connected DG algebra. If the cohomology
algebra $H(A)$ is a left AS-Gorenstein algebra, then $A$ is a left
AS-Gorenstein DG algebra.
\end{prop}

\proof Consider the Eilenberg-Moore spectral sequence (\cite{KM})
$$E_2^{p,q} = \Ext^p_{H(A)}(k,H(A))^q\Longrightarrow \Ext_A^{p+q}(k,A)$$
where the index $p$ in $\Ext^p_{H(A)}(k,H(A))^q$ is the usual
homological degree and $q$ is the grading induced from the gradings
of ${}_Ak$ and $H(A).$ If the cohomology spectral sequence is
regular, then it is complete convergent (\cite{We}). If $H(A)$ is
AS-Gorenstein, then by definition there exist some integers $d$ and
$l$ such that
$$
\Ext_{H(A)}^n(k, H(A)) =
\begin{cases}
0 & n \neq d \\
k[l] &  n=d.
\end{cases}
$$
Then it is routine to see that
$$\Ext^n_A(k,A) =
\begin{cases}
0 & n \neq d+l \\
k &  n=d+l.
\end{cases}
$$
Hence $\RHom_A(k,A) \cong s^nk$ for $n=d+l$. $\square$

We don't know whether the converse of Proposition \ref{prop9} is
true or not. If $A$ is a connected graded algebra, viewed as a DG
algebra with trivial differential, then $A$ is an AS-Gorenstein DG
algebra if and only if $A$ satisfies AS-Gorenstein condition in the
usual sense. The AS-Gorenstein property is invariant under
quasi-isomorphism.
\begin{prop}Let $f:A\to A'$ be a quasi-isomorphism of connected DG
algebras. Then $A$ is left AS-Gorenstein (AS-regular) if and only if
$A'$ is.
\end{prop}\proof The proof is similar to that of Proposition
\ref{prop4}. $\square$

\begin{lem}\label{lemm} Let $A$ be an AS-regular DG algebra. Suppose
$\RHom_{A\!^{op}}(k,A) \cong s^lk$ for some integer $l$. Let $P \to
k_A$ be a minimal semifree resolution of $k_A$ with a semifree
filtration
$$0\subseteq P(0)\subseteq P(1)\subseteq\cdots\subseteq P(n)$$ such
that $P(n)=P$ and $P(n)/P(n-1)\neq0$. Then $P(n)/P(n-1)=A[-l]$.
\end{lem}

\proof There are finite dimensional graded vector spaces $V(0),\
V(1),\cdots,V(n)$ such that $P(i)/P(i-1)=V(i)\otimes A$ for all $0
\leq i\leq n$ ($P(-1)=0$). As graded $A$-modules
$P=\bigoplus_{i=0}^n V(i)\otimes A$. Hence
$\Hom_A(P,A)=\bigoplus_{i=0}^nA\otimes V(i)^\#$ as graded left
$A$-modules. Let $\{x_1,\cdots,x_t\}$ be a homogeneous basis of
$V(n)$. Let $d$ be the differential of $\Hom_A(P,A)$ induced by the
differentials of $P$ and $A$. For any $1\leq s\leq t$, define a
graded right $A$-module morphism
$$f_s:P=\oplus_{i=0}^n V(i)\otimes A\longrightarrow A$$ by sending $x_s$
to the identity of $A$, $x_j$ to zero for $j \neq s$, and sending
$V(r)$ to zero for all $r<n$. One can see that $f_1,\cdots,f_t$ so
defined are cocycles of the cochain complex $\Hom_A(P,A)$. Since $P$
is minimal, $d(g)(x_j)=(d_A \,g -(-1)^{|g|}g d_P)(x_j)\in A^{\ge1}$
for any homogeneous element $g \in \Hom_A(P,A)$ and $x_j$. Hence any
$f_j (1 \leq j \leq t)$ can not be a coboundary. By hypothesis
$\RHom_{A\!^{op}}(k,A)\cong s^lk$, which forces $\dim V(n)=1$ and
the degree of non-zero elements in $V(n)$ is $-l$, that is,
$P(n)/P(n-1)\cong s^{-l}A$. $\square$

The following proposition is a special case of \cite[Theorem
9.8]{LP2}.

\begin{prop}\label{prop7} Let $A$ be a Koszul DG algebra with
Ext-algebra $E=\Ext_A^*(k,k)$. Then $A$ is right AS-regular if and
only if $E$ is Frobenius.
\end{prop}
\proof Suppose that $A$ is right AS-regular. Then $k_A\in
\mathcal{D}^c(A\!^{op})$, which is equivalent to ${}_Ak\in
\mathcal{D}^c(A)$. Hence $E$ is finite dimensional. Since $A$ is
Koszul, $\RHom_{A\!^{op}}(k,A)\cong k$ by Lemma \ref{lemm}. By
Theorem \ref{thm6},
$$\begin{array}{ccl}
    \Ext^n_{E\!^{op}}(k,E)&=&\Hom_{\mathcal{D}^b(\text{mod-}E\!^{op})}(k,E[-n])\\
 &\cong & \Hom_{\mathcal{D}^c(A\!^{op})}(\mathcal{F}(E),\mathcal{F}(k)[-n]) \\
     &\cong& \Hom_{\mathcal{D}^c(A\!^{op})}(k,A[-n]) \\
  & = &  \Ext_{A\!^{op}}^{n}(k,A).
  \end{array}
$$ Hence $\Ext^n_{E\!^{op}}(k,E)=0$ for $n\neq0$. Let
$$0\longrightarrow E_E\longrightarrow I^0\longrightarrow
I^1\longrightarrow\cdots$$ be a minimal injective resolution of
$E_E$. Since $E$ is finite dimensional and local, all the injective
modules $I^n$'s are finite dimensional. Hence
$0=\Ext^n_{E\!^{op}}(k,A)\cong socI^n$ for all $n\ge1$, and
$\Hom_{E\!^{op}}(k,E)=k$. We get $I^n=0$ for all $n\ge1$ and
$I^0=E^*$. Therefore we have a right $E$-module isomorphism $E\cong
E^*$, that is, $E$ is a Frobenius algebra.

Conversely, if $E$ is Frobenius, then it is finite dimensional, and
hence $k_A\in \mathcal{D}^c(A^{op})$. Since $E$ itself is injective
and local, it follows $\Ext^n_{E\!^{op}}(k,E)=0$ for $n\ge1$ and
$\Ext^0_{E\!^{op}}(k,E)=k$. Hence $\Ext_{A\!^{op}}^{n}(k,A)\cong
\Ext^n_{E\!^{op}}(k,E)=0$ for $n\neq0$ and
$\Ext^0_{A\!^{op}}(k,A)\cong k$. Then $\RHom_{A\!^{op}}(k,A)\cong
k$, and hence $A$ is AS-regular. $\square$

In \cite[Theorem 9.8]{LP2}, a  more general case of the above
proposition is proved with some locally finite conditions.

\begin{cor}\label{cor6}Let $A$ be a Koszul DG algebra. Then $A$ is right
AS-regular if and only if $A$ is left AS-regular.
\end{cor}
\proof Note that $E^{op} \cong \Ext_{A^{op}}^*(k,k).$ $\quad
\square$

Next we are going to deduce a result similar to the classical BGG
correspondence.

\begin{lem}\label{lem8} Let $A$ be a connected DG algebra such that
$k_A\in \mathcal{D}^c(A\!^{op})$. Then the full triangulated
subcategory $\langle k_A\rangle$ of $\mathcal{D}^c(A\!^{op})$
generated by $k_A$, is equal to $\mathcal{D}_{fd}(A\!^{op})$, the
full subcategory of $\mathcal{D}^c(A\!^{op})$ consisting of DG
modules $M$ such that $\dim H(M)<\infty$.
\end{lem}

\proof For any DG module $M$, temporarily we write
$$\ell(M)=\sup\{i\,|\,H^i(M)\neq0\}-\inf\{i\,|\,H^i(M)\neq0\}$$ and
$$\lambda(M)=\sup\{i\,|\,M^i\neq0\}-\inf\{i\,|\,M^i\neq0\}.$$
We prove the lemma by an induction on $\ell(M)$.
Let $M$ be a DG $A$-module with $\dim H(M)<\infty$. Without loss of
generality, we may assume that $H^i(M)=0$ for $i<0$ or $i>n$ for
$n=\ell(M)$. Since $A$ is connected, by suitable truncations, we may
assume that $M$ is concentrated in degrees $0 \leq i \leq n.$
If $\ell(M)=0$, then $M$ is isomorphic in $\mathcal{D}^c(A\!^{op})$
to a DG module $N$ with $\lambda(N)=0$, which is a direct sum of
finite copies of $k_A$, and hence is in $\langle k_A\rangle$. Now
suppose that each DG module $M$ with $\dim H(M)<\infty$ and
$\ell(M)< n$ is in $\langle k_A\rangle$. If $M$ is a DG module with
$\dim H(M)<\infty$ and $\ell(M)=n$, without loss of generality, we
may assume $M^0\neq0$ and $M^n\neq0$, and $M^i=0$ for $i<0$ and
$i>n$. Then the vector space $M^n$ has a decomposition
$M^n=d(M^{n-1})\op K$ for some subspace $K$. Since $\dim
H(M)<\infty$, then $\dim K<\infty$. Taking $K$ as a DG $A$-module
concentrated on degree zero, then we have an exact sequence of DG
modules
$$0\longrightarrow s^nK\longrightarrow M\longrightarrow
M/s^nK\longrightarrow0.$$ Now $\ell(s^nK)=0$ and $\ell(M/s^nK)\leq
n-1$. By the induction hypothesis, both $s^nK$ and $M/s^nK$ are
objects in $\langle k_A\rangle$, and hence $M$ is in $\langle
k_A\rangle$. Therefore $\langle k_A\rangle =
\mathcal{D}_{fd}(A\!^{op})$. $\square$

\begin{thm}[BGG Correspondence]\label{thm7} Let $A$ be a Koszul DG
AS-regular algebra with Ext-algebra $E=\Ext_A^*(k,k)$. Then there is
a duality of triangulated categories
$$\overline{\text{\rm mod-}}\,E\!^{op} \rightleftarrows \mathcal{D}^c(A\!^{op}) / \mathcal{D}_{fd}(A\!^{op}).$$
\end{thm}

\proof By the Koszul Duality (Theorem \ref{thm6}), There is  a
duality of triangulated categories
$$\mathcal{D}^b(\text{\rm
mod-}\,E\!^{op}) \rightleftarrows \mathcal{D}^c(A\!^{op});$$ and
under this duality the object $E_E \in \mathcal{D}^b(\text{\rm
mod}E\!^{op})$ is corresponding to the object $k_A \in
\mathcal{D}^c(A\!^{op})$ by (\ref{eq4}). Hence there is a duality
$$\mathcal{D}^b(\text{\rm mod-}E\!^{op})/\langle E_E\rangle \rightleftarrows
\mathcal{D}^c(A\!^{op})/\langle k_A\rangle,$$ where $\langle
E_E\rangle$ is the full triangulated subcategory of
$\mathcal{D}^b(\text{\rm mod-}E\!^{op})$ generated by $E_E$. Since
$E$ is a finite dimensional local algebra with $E/J(E) \cong k$
(Theorem \ref{thm2}), all finitely generated projective $E$-modules
are free. Therefore  $\langle E_E\rangle =
\mathcal{D}^b(\text{proj}\,E\!^{op})$, where $\text{proj}\,E\!^{op}$
is the category of all finitely generated right projective
$E$-modules. Hence
$$\mathcal{D}^b(\text{\rm mod-}E\!^{op})/\langle E_E\rangle = \mathcal{D}^b(\text{\rm
mod-}E\!^{op})/\mathcal{D}^b(\text{proj}\, E\!^{op}).$$ By
Proposition \ref{prop7}, $E$ is Frobenius, and hence (\cite{Bel})
$$\mathcal{D}^b(\text{\rm mod-}E\!^{op})/\mathcal{D}^b(\text{proj}\,
E\!^{op})\cong\overline{\text{\rm mod-}}E\!^{op}.$$ On the other
hand, by Lemma \ref{lem8}
$$\mathcal{D}^c(A\!^{op})/\langle
k_A\rangle = \mathcal{D}^c(A\!^{op})/\mathcal{D}_{fd}(A).$$ In
summary, there is an a duality of triangulated categories
$$\overline{\text{\rm mod-}}E^{op} \rightleftarrows
\mathcal{D}^c(A\!^{op})/\mathcal{D}_{fd}(A\!^{op}).\quad \square$$

Since $E$ is finite dimensional, there is an equivalence form of the
BGG Correspondence.

\begin{thm}\label{thm10} Let $A$ be a Koszul DG
AS-regular algebra with Ext-algebra $E=\Ext_A^*(k,k)$. Then there is
an equivalence of triangulated categories
$$\overline{\text{\rm mod-}}\,E \cong \mathcal{D}^c(A\!^{op}) / \mathcal{D}_{fd}(A\!^{op}).$$
\end{thm}

\section{BGG correspondence on Adams connected DG algebras}

Many examples of the DG algebra from algebraic geometry and
algebraic topology admit an extra grading. Let $A=\op_{i,j\in\Bbb
Z}A^i_j$ be a bigraded space.  An element $a\in A^i_j$ is of degree
$(i,j)$. The second grading is usually called {\it Adams grading}
(\cite{KM}, \cite{LPWZ}). A  DG algebra $(A,d)$ is called a DG
algebra with Adams grading if  $A$ is bigraded and the differential
$d$ is of degree $(1,0)$ (i.e., $d$ preserves Adams grading). A DG
module over a DG algebra with Adams grading is bigraded and the
differential preserves the second grading. A DG algebra $A$ with
Adams grading is {\it augmented} if there is an augmentation map
$\varepsilon:A\to k$ of degree $(0,0)$. A DG algebra with Adams
grading $A$ is said to be {\it Adams connected} if (1) $A^i_j=0$ for
$i<0$ or $j<0$, and (2) $A^0_0=k$, $A^0_j=0$ and $A^i_0=0$ for
$i,j\neq0$. All Adams connected DG algebras are augmented.

Similarly, we define coaugmented DG coalgebras with Adams grading.

In this section, all the DG algebras and DG coalgebras involved are
with Adams grading. For simplicity, we call a DG algebra (coalgebra)
with Adams grading an {\it Adams DG algebra (coalgebra)}.

It is not hard to see that the bar (cobar) construction (see Section
1) of an (a) (co)augmented Adams DG algebra (coalgebra) is an Adams
DG coalgebra (algebra). The canonical twisting cochain (see Section
4) $\tau_0$ from a cocomplete Adams DG coalgebra $C$ to $\Omega(C)$
is of degree $(1,0)$.

Let $A$ be an Adams DG algebra, and let $\mathcal{AC}_{dg}(A)$
$(\mathcal{AC}_{dg}(A\!^{op}))$ be the category of left (right) DG
$A$-modules with morphisms of degree $(0,0)$. We use $s^i=[i]$ to
denote the $i$-th shift functor on the first grading and use
$s^{-j}=(j)$ to denote the $j$-th shift functor on the Adams
grading. Let $\mathcal{AD}_{dg}(A)$ be the derived category of
$\mathcal{AC}_{dg}(A)$. Denote $\mathcal{AD}^c(A)$
($\mathcal{AD}^c(A\!^{op})$) as the full triangulated subcategory of
$\mathcal{AD}_{dg}(A)$ ($\mathcal{AD}_{dg}(A\!^{op})$) generated by
${}_AA$ ($A_A$). Let $M$ and $N$ be objects in
$\mathcal{AC}_{dg}(A)$, we use
$$\mathcal{E}xt_A^{i,j}(M,N)=\Hom_{\mathcal{AD}_{dg}(A)}(M,N[-i](j))$$
to denote the derived functor. Then
$$\mathcal{E}xt_A^{*,*}(M,N)=\bigoplus_{i,j\in\Bbb Z}\mathcal{E}xt_A^{i,j}(M,N)$$
is a bigraded space. In particular, if $A$ is an augmented Adams DG
algebra, then $\mathcal{E}=\mathcal{E}xt_A^{*,*}(k,k)$ is a bigraded
algebra. For convenience, we usually write
$\mathcal{E}^i_j=\mathcal{E}xt_A^{i,j}(k,k)$.

The results obtained in previous sections can be easily generalized
to Adams DG algebras. Hence in this section, we only state the
results without giving proofs. More general results can be found in
\cite[Section 10]{LP2}, with some locally finite conditions.

Let $A$ be an Adams connected DG algebra and $M$ be a bounded below
DG module over $A$. Then there is a minimal semifree resolution (the
construction is similar to \cite[Theorem IV.3.7]{KM}) $P\to M$ in
$\mathcal{AC}_{dg}(A)$ (see also \cite{MW}).

\begin{defn} {\rm Let $A$ be an Adams connected DG algebra. It is
called a {Koszul} Adams DG algebra if
$\mathcal{E}xt_A^{i,*}(k,k)=\bigoplus_{j \in \Bbb
Z}\mathcal{E}xt_A^{i,j}(k,k)=0$ for all $i\neq0$.}
\end{defn}

It is not hard to see that, if $A$ is a Koszul Adams DG algebra,
then its Ext-algebra $\mathcal{E}=\mathcal{E}xt_A^{*,*}(k,k)$ has
the property that $\mathcal{E}^i_j=0$ for $i\neq0$ or $j>0$. Hence
$\mathcal{E}$ is a negatively graded algebra. Comparing with Theorem
\ref{thm1}, we have the following.

\begin{prop}\label{10} Let $A$ be a Koszul Adams DG algebra, and let
$S_j=\mathcal{E}^0_{-j}$. Then $S=\op_{j\ge0}S_j$ is a connected
graded algebra. If in addition ${}_Ak\in \mathcal{AD}^c_{dg}(A)$,
then $S$ is a finite dimensional graded algebra.
\end{prop}

We also have the following form of Lef\`{e}vre-Hasegawa's Theorem.

\begin{thm} Let $C$ be a cocomplete Adams DG coalgebra, $B$ an augmented Adams DG
algebra and $\tau:C\to B$ is a twisting cochain of degree $(1,0)$.
The following are equivalent
\begin{itemize}
    \item [(i)] The map $\tau$ induces a quasi-isomorphism $\Omega(C)\to
    B$;
    \item [(ii)] The adjunction map $$B\otimes_\tau C\otimes_\tau B\to B$$ is a
    quasi-isomorphism;
    \item [(iii)] There is an equivalence of triangulated categories
     $$\mathcal{AD}_{dg}(C)\rightleftarrows\mathcal{AD}_{dg}(B\!^{op})$$ where $\mathcal{AD}_{dg}(C)$ is the coderived
     category over the cocomplete Adams DG algebra $C$. $\square$
\end{itemize}
\end{thm}

If ${}_Ak\in\mathcal{AD}^c(A)$, then $\mathcal{E}$ is finite
dimensional, and hence the graded vector space dual $\mathcal{E}^\#$
is finite dimensional coalgebra. By applying the above theorem and
notice that a DG comodule over the Adams DG coalgebra
$\mathcal{E}^\#$ is exactly a complex of graded comodules over
$\mathcal{E}^\#$, we have the following proposition which is
analogous to Theorem \ref{thm6}.

\begin{prop}\label{prop8} Let $A$ be a Koszul Adams DG algebra. If
${}_Ak\in\mathcal{AD}^c(A)$, then there is a duality of triangulated
categories
$$\mathcal{D}^b(\text{\rm
grmod-}\mathcal{E}^{op}) \darrow {\mathcal{F}} {\mathcal{G}}
\mathcal{AD}^c(A\!^{op}). \quad\square$$
\end{prop}

It is convenient for us to deal with the positively graded algebra
$S$, rather than the negatively graded algebra $\mathcal{E}$. The
abelian category $\text{grmod-}\,\mathcal{E}^{op}$ is equivalent to
$\text{grmod-}S^{op}$ of finitely generated right $S$-modules. We
have the following Koszul duality theorem of Adams DG algebras.

\begin{thm}\label{thm8} Let $A$ be a Koszul Adams DG algebra. Let $S$ be the graded algebra such that
$S_j=\mathcal{E}xt^{0,-j}_A(k,k)$. If ${}_Ak \in \mathcal{AD}^c(A)$,
then there is a duality of triangulated categories
$$\mathcal{D}^b(\text{\rm
grmod-}S^{op}) \darrow {\psi} {\phi} \mathcal{AD}^c(A\!^{op}).
\quad\square$$
\end{thm}

To establish a version of the BGG correspondence, we need the
concept of AS-Gorenstein Adams DG algebra which is first introduced
in \cite{LPWZ}.

\begin{defn} {\rm Let $A$ be an Adams connected DG algebra. It is called
an {\it AS-Gorenstein} Adams DG algebra if
$\text{RHom}_{A\!^{op}}(k,A)\cong k[r](s)$. Moreover if
$k_A\in\mathcal{AD}^c(A\!^{op})$, then $A$ is called an {\it
AS-regular} Adams DG algebra.}
\end{defn}

The following proposition is proved in \cite{LPWZ} by using
$A_\infty$-algebra. Also one can give a proof by using Theorem
\ref{thm8}.

\begin{prop}Let $A$ be a Koszul Adams DG algebra. Then $A$ is
AS-regular if and only if its Ext-algebra $\mathcal{E}$ is
Frobenius.
\end{prop}

Now we can state the BGG correspondence on Adams DG algebras.

\begin{thm}\label{thm9} Let $A$ be a Koszul
AS-regular Adams DG algebra and $S$ be the connected graded such
that $S_j=\mathcal{E}xt_A^{0,-j}(k,k)$. Then there is an duality of
triangulated categories
\begin{equation}\label{eq6}
\overline{\text{\rm grmod-}}S^{op}\cong
\mathcal{AD}^c(A\!^{op})/\mathcal{AD}_{fd}(A\!^{op}),
\end{equation}
 where $\mathcal{AD}_{fd}(A\!^{op})$ is the full triangulated
subcategory of $\mathcal{AD}^c(A\!^{op})$ consisting of objects $M$
such that $\dim H(M)<\infty$.
\end{thm}

Now let $R$ be a noetherian connected graded algebra. Let $A$ be the
Adams connected DG algebra with trivial differential by taking
$A^i_i=R_i$ and $A^i_j=0$ if $i\neq j$. If $R$ is a Koszul algebra,
then it is not hard to see that $A$ is a Koszul Adams DG algebra.
Moreover, $Ext_A^{0,-j}(k,k)= R^!_j$ for all $j\ge0$, i.e.,
$S=R^!=E(R)=\Ext_R^*(k,k)$. Suppose that $\text{gl}.\dim R<\infty$.
Then $S=E(R)$ is finite dimensional. Since $E(R)$ is finite
dimensional $\text{\rm grmod-}E(R)^{op}$ is dual to $\text{\rm
grmod-}E(R)$. Hence $\mathcal{D}^b(\text{\rm
grmod-}S^{op})=\mathcal{D}^b(\text{\rm grmod-}E(R)^{op})$ is dual to
$\mathcal{D}^b(\text{\rm grmod-}E(R))$. Let us inspect the category
$\mathcal{AD}^c(A\!^{op})$ in Theorem \ref{thm8}. Since the
differential of $A$ is trivial and $A$ is concentrated in the
diagonal of the first quadrant, the triangulated category
$\mathcal{AD}_{dg}(A\!^{op})$ is naturally equivalent to the derived
category $\mathcal{D}(\text{Grmod-}R^{op})$ of the category
$\text{Grmod-}R$ of right graded $R$-modules. Under this
equivalence, $A_A$ is corresponding to $R_R$ in
$\mathcal{D}(\text{Grmod-}R^{op})$. Hence $\mathcal{AD}^c(A\!^{op})$
is equivalent to the full triangulated subcategory of
$\mathcal{D}(\text{Grmod-}R^{op})$ generated by $R_R$ (closed under
the shifts on the grading of $R_R$), which is equivalent to
$\mathcal{D}^b(\text{proj}\,R^{op}),$ the bounded derived category
of finitely generated  graded projective right $R$-modules. Since
$R$ is noetherian and has finite global dimension,
$\mathcal{D}^b(\text{proj}R^{op})$ is equivalent to
$\mathcal{D}^b(\text{grmod-}R^{op}),$ the bounded derived category
of finitely generated  graded right $R$-modules. In summary we have
the equivalence (which is established in \cite{BGS}) of triangulated
categories if $R$ is noetherian and of finite global dimension
\begin{equation*}
    \mathcal{D}^b(\text{\rm grmod-}\,E(R))\cong
\mathcal{D}^b(\text{grmod-}\,R^{op}).
\end{equation*}

Moreover, we assume that $R$ is a noetherian Koszul AS-regular
algebra. Then the Adams connected DG algebra $A$ is Koszul Adams
AS-regular DG algebra. Hence in the left hand of (\ref{eq6}),
$\overline{\text{\rm grmod-}}S^{op}$ is dual to $\overline{\text{\rm
grmod-}}\,E(R)$. Since $\mathcal{AD}^c(A\!^{op})$ is equivalent to
$\mathcal{D}^b(\text{grmod-}R^{op})$, the full triangulated
subcategory $\mathcal{AD}_{fd}(A\!^{op})$ is equivalent to
$\mathcal{D}^b_{fd}(\text{grmod-}R^{op})$, the triangulated
subcategory consisting of objects $X$ such that $HX$ is finite
dimensional. Hence in the right hand of (\ref{eq6}),
$$\mathcal{AD}^c(A\!^{op})/\mathcal{AD}_{fd}(A\!^{op})\cong
\mathcal{D}^b(\text{grmod-}R^{op})/\mathcal{D}^b_{fd}(\text{grmod-}R^{op})$$
which is equivalent to $\mathcal{D}^b(\text{qgr}R^{op})$ by
\cite{Miy}. In summary we get the BGG correspondence established in
\cite{MS}
$$\overline{\text{\rm grmod-}}\,E(R)\cong \mathcal{D}^b(\text{qgr}R^{op}).$$

\vspace{2mm}

\section*{Acknowledgments}

This research is supported by the NSFC (key project 10331030) and
supported by the Cultivation Fund of the Key Scientific and
Technical Innovation Project (No. 704004), Doctorate Foundation (No.
20060246003), Ministry of Education of China.

\bibliography{}

\end{document}